\documentclass[12pt, reqno]{amsart}

\usepackage[hmargin=0.8in,height=8.6in]{geometry}
\usepackage{amssymb,amsthm, times}
\usepackage{delarray,verbatim}
\usepackage{ifpdf}
\ifpdf
\usepackage[pdftex]{graphicx}
 \else
\usepackage[dvips]{graphicx}
 \fi

\usepackage{bm}

\usepackage{ifpdf}
\usepackage{color}
\definecolor{webgreen}{rgb}{0,.5,0}
\definecolor{webbrown}{rgb}{.8,0,0}
\definecolor{emphcolor}{rgb}{0.95,0.95,0.95}

\usepackage[pagebackref]{hyperref}
\usepackage{hyperref}
\hypersetup{%
          colorlinks=true,
          linkcolor=webbrown,
          filecolor=webbrown,
          citecolor=webgreen,
          breaklinks=true}
\ifpdf \hypersetup{pdftex,
            pdfstartview=FitH, 
            bookmarksopen=true,
            bookmarksnumbered=true
} \else \hypersetup{dvips} \fi

\newcommand{\lapinv}{\Phi(q)}

\numberwithin{equation}{section}

\newtheorem{theorem}{Theorem}[section]
\newtheorem{proposition}{Proposition}[section]

\newtheorem{remark}{Remark}[section]
\newtheorem{lemma}{Lemma}[section]

\newtheorem{assump}{Assumption}[section]
\newtheorem{definition}{Definition}[section]
\numberwithin{remark}{section} \numberwithin{proposition}{section}
\numberwithin{corollary}{section}
\newcommand {\R}{\mathbb{R}}

\newcommand {\F}{\mathcal{F}}

\newcommand {\p}{\mathbb{P}}

\newcommand {\E}{\mathbb{E}}

\newcommand{\diff}{{\rm d}}

\newcommand{\lev}{L\'{e}vy }

\title[Inventory Control for Spectrally Positive L\'evy Demand Processes]{Inventory Control for Spectrally Positive L\'evy Demand Processes}
\thanks{This version: \today. }
\thanks{$*$\, (corresponding author) Department of Mathematics,
Faculty of Engineering Science, Kansai University, Suita-shi, Osaka 564-8680, Japan. Email: \mbox{{\em
kyamazak@kansai-u.ac.jp}} Phone: +81-6-6368-1527.  }
\author[K. Yamazaki]{Kazutoshi Yamazaki$^*$}
\date{}
\begin{document}

\begin{abstract}
A new approach to solve the continuous-time stochastic inventory problem using the fluctuation theory of \lev processes is developed.  This approach involves the recent developments of the scale function that is capable of expressing many fluctuation identities of spectrally one-sided \lev processes.  For the case with a fixed cost and a general spectrally positive \lev demand process, we show the optimality of an $(s,S)$-policy.  The optimal policy and the value function are concisely expressed via the scale function.  Numerical examples under a \lev process in the $\beta$-family with jumps of infinite activity are provided to confirm the analytical results.
Furthermore, the case with no fixed ordering costs is studied.
%
%
\\
\noindent \small{\noindent  AMS 2010 Subject Classifications: 60G51, 93E20, 49J40 
}\\
\end{abstract}

\maketitle

\section{Introduction}

In this study, we introduce the \emph{fluctuation theory} of \lev processes to solve inventory control problems.  In the continuous-time model, the majority of existing studies use a certain type of \lev process as a demand process; typically it is modeled by Brownian motion,  a compound Poisson process,  or a mixture of the two. Pursuing the connections between the inventory theory and the L\'evy-process theory is naturally of great importance. However, 
recent advances on \lev processes have not been used to study inventory models.
Therefore, our objective in the current study is to fill this void. Toward this end, we introduce the so-called \emph{scale function} to inventory control, which plays a key role in the fluctuation theory of \lev processes.   We demonstrate both analytically and numerically that it is in fact a powerful tool to solve inventory control problems.

The common aim of inventory control is to derive the \emph{optimal replenishment policy} to minimize both  inventory and ordering costs.  This study focuses on the discounted continuous-time model with fixed and proportional costs.  We assume that the order quantity is continuous and back-orders are allowed. Furthermore, we assume the absence of lead time, perishability, and lost sales.  This setting is the same as that in the seminal study by \cite{Bensoussan_2009} where they solved (complementing \cite{Bensoussan_2005}) the case where the demand arrives as a combination of Brownian motion and a compound Poisson process.  In the current study, a scale-function-based approach is used to solve the problem for any choice of \emph{spectrally positive \lev demand process} i.e.\ a \lev process with only upward jumps. 

To our best knowledge, this is the first study on inventory models where the theory of scale functions is applied to derive the optimal solution.  We focus on the aforementioned simple setting and demonstrate that the existing known properties of the scale function can be efficiently used to follow the classical  \emph{guess and verify} procedure in a straightforward fashion.
Of course, applying the same technique in other inventory models can be of great importance. We expect that this study can potentially serve as a guide on how to tackle these open problems using the theory of scale functions.
For a detailed review of inventory control problems, we refer the readers to \cite{Bensoussan_2005,Bensoussan_Tapiero_1982,Presman_Sethi_2006} and the recent book by Bensoussan \cite{Bensoussan_book_2011}.

Over the last decade or two, significant developments in the fluctuation theory of \lev processes have been presented (see, e.g., the textbooks by Bertoin \cite{Bertoin_1996}, Doney \cite{Doney_2007}, and Kyprianou \cite{Kyprianou_2006}).   Whereas 
the majority of \lev processes still
 remain to be analytically tractable,  numerous computations have been made possible when it comes to the spectrally one-sided \lev process (a \lev process with only upward or downward jumps) because of the development of the scale function.  

For every spectrally one-sided \lev process, the corresponding scale function, which is concisely defined in terms of the Laplace exponent, plays a great role. Without specifying a particular process type, it can efficiently express, e.g., hitting time probabilities, resolvent (potential) measures, and overshoot/undershoot distributions.  It has gradually come into use in several stochastic models. In particular, it is now considered as a key tool in insurance mathematics; the classical formulation that models the surplus of an insurance company using a compound Poisson process is now being replaced by the spectrally negative \lev model. Other stochastic models where the scale function plays a key role include, e.g., optimal stopping \cite{Egami-Yamazaki-2010-1, Kyprianou_Ott_2014, leung_yamazaki_2013, Ott_2013, Yamazaki_2014_contraction}, stochastic games    \cite{Baurdoux2008, BauKyrianouPardo2011, Leung_Yamazaki_2011, Hernandez_Yamazaki_2013}, and mathematical finance  \cite{Kyprianou_Surya_2007, surya_yamazaki_2014}.

In this study, we begin with the same formulation as in \cite{Bensoussan_2005} and solve for the general spectrally positive \lev process.  We have no difficulty in conjecturing that the form of an optimal policy is again of the $(s,S)$-type, i.e., bringing the inventory level up to $S$ whenever it goes below $s$ is optimal.  However, the scale-function-based approach is worth studying for various reasons.  Here we summarize its advantages over the classical approach that involves solutions to integro-differential equations (IDEs).

\underline{Generalization of the demand process:} 
Essentially, all existing L\'evy-based inventory models employ the \lev processes with upward jumps of finite activity [as defined in \eqref{compound_poisson}], i.e.,  the demand arrives as a (drifted) compound Poisson process with or without a Brownian motion. The set of these \lev processes excludes those with infinite activity/variation. In other words, the existing models do not cover  well-known processes such as the (spectrally positive versions of) variance gamma, CGMY, and normal inverse Gaussian processes as well as classical ones as the gamma process and a subset of stable processes.  For these processes, the IDE involves integration with respect to infinite measures, and hence, the IDE-based approach can become intractable.
On the other hand, the approach that uses the scale function can accommodate these processes without  any additional work. 

 \lev processes of infinite activity/variation have been applied in the field of finance and insurance.    Empirical evidence shows that asset prices can be more precisely modeled by jumps of infinite activity as opposed to being modeled by Brownian motion (see the introduction in \cite{CGMY_2002}).  Because the demand is closely linked to the price of an item, modeling inventory systems using these processes also makes sense.  
 This reason motivates us to work on the inventory control problem for a general spectrally positive \lev demand process.



\underline{Conciseness of arguments:}  The scale function can explicitly express the expected total costs under the $(s,S)$-policy.  Using this function, the optimal solution can be derived in a straightforward manner.  As typically conducted in solving stochastic control problems, the optimal solution is first ``guessed," and its optimality is then confirmed by ``verification" arguments. The known properties of the scale function help one achieve these objectives.

The first guessing step is reduced to choosing the pair $(s^*, S^*)$ by imposing  the following two conditions on the expected cost function as follows:
\begin{enumerate}
\item The continuous (resp.\ smooth) fit condition at the lower threshold $s$ when the demand process is of  bounded (resp.\ unbounded) variation, 
\item The condition at the upper threshold $S$ such that the slope at $S$ is equal to the negative of the unit proportional cost.
\end{enumerate}
Because of the (semi-)analytical expression of the expected costs under the $(s,S)$-policy, these two conditions can be concisely rewritten by the scale function. For the results in this study,  $(s^*,S^*)$ become the zeros of the function $\mathcal{G}(s,S)$  defined in \eqref{def_G} and that of its derivative with respect to the second argument.   Using this graphical interpretation and taking advantage again of the certain known properties of the scale function, the existence of such  $(s^*, S^*)$ can be confirmed.

The second verification step reduces to indicating that the candidate value function is sufficiently smooth and  satisfies the \emph{quasi-variational inequalities (QVIs)}.  The former can be easily confirmed because of the known smoothness/continuity properties of the scale function that depend on the path variation of the process.  Unfortunately,  part of the latter can become difficult to verify. However, in this study, the optimality in fact holds because of several known facts on the scale function and the results obtained in \cite{Bensoussan_2009}.

\underline{Computability:}  The derived optimal thresholds $(s^*, S^*)$ and the associated optimal value function are concisely written via the scale function.   Hence, the computation of these factors is essentially equivalent to that of the scale function.  Because the scale function is defined by its Laplace transform written in terms of the Laplace exponent,  the Laplace transform needs to be inverted either analytically or numerically to compute it.  

Fortunately, some important classes of \lev processes have rational forms of Laplace exponents. For these processes, analytical forms of scale functions can be easily obtained by partial fraction decomposition. In Section \ref{subsection_examples_scale_function}, we show examples of such processes. Among others, the phase-type \lev process of \cite{Asmussen_2004} is of great interest both in analytical and numerical aspects. It admits a rational form of Laplace exponent and is known to be \emph{dense} in the set of all \lev processes.  This means that, it has an explicit form of scale function, and, more importantly,  any scale function can be approximated by the scale function of this form. Egami and Yamazaki \cite{Egami_Yamazaki_2010_2} conducted a sequence of numerical experiments to confirm the accuracy of this approximation.


Alternatively, the scale function can always be directly computed via numerical Laplace inversion.  As discussed in Kuznetsov et al.\ \cite{Kuznetsov2013}, a scale function can be written as the difference between an exponential function (whose parameter is defined by $\Phi(q)$ in the current paper) and the resolvent (potential) term [see the third equation in \eqref{resolvent_density} below].  Hence, the computation is reduced to that of the resolvent term. It is a bounded function that asymptotically converges to  zero, and hence, numerical Laplace inversion can be quickly and accurately conducted. For more details, we refer the readers to Section 5 of \cite{Kuznetsov2013}.

To confirm the analytical results obtained in this study, we provide numerical examples with a quadratic inventory cost and a demand process in the $\beta$-family in \cite{Kuznetsov_2010_2} with jumps of infinite activity.  We can see that the  optimal levels $(s^*,S^*)$ and the value function can be instantaneously computed.  

\underline{Sensitivity analysis:} 
Because of the analytical expressions of optimal thresholds and the value function,  sensitivity analysis can be more directly conducted.  

The sensitivity with respect to the parameter of the underlying \lev process is equivalent to that of the scale function.  For example, the smoothness of the optimal value function is directly linked to the asymptotic behavior of the scale function near zero.  In our numerical results, we consider the cases where the diffusion coefficient is zero and positive, and analyze how the optimal solutions differ. 

Furthermore, the sensitivity with respect to the fixed and unit proportional costs is of great interest.  In particular, the forms of the optimal solutions are different depending on the existence of a fixed ordering cost. The distance between the optimal thresholds $s^*$ and $S^*$ is expected to shrink as the fixed cost decreases.
 In this study, we also consider the case with no fixed ordering cost and show that the optimal policy is of the barrier type.  Using the fluctuation theory of reflected \lev processes as in \cite{Avram_et_al_2007} and \cite{Pistorius_2004}, the value function can again be written using the scale function.  We  numerically confirm that, as the fixed cost decreases to zero, the optimal $(s^*,S^*)$-policy converges to the optimal barrier strategy. 


\vspace{0.3cm}

Before closing the introduction, we briefly review the scale-function approach used in other stochastic control problems and discuss the similarities and differences with the results of this study.  

The most relevant factor is the optimal dividend problem (de Finetti's problem) in insurance mathematics. After the development of the  fluctuation theory of reflected \lev processes, many authors have applied the results in the spectrally negative \lev model of de Finetti's problems \cite{Kyprianou_Palmowski_2007,Loeffen_2008, Loeffen_2009, Loeffen_2009_2} (see also \cite{Bayraktar_2012,Bayraktar_2013, Yin_Wen_2013} for the spectrally positive \lev models).  The main difference with our inventory control problem is that whereas our problem has an infinite time horizon,  de Finetti's problem is terminated at the ruin time (or the first time the controlled process goes below the zero level).  In our problem, the absence of a ruin makes the problem both easier and more difficult.  Whereas in de Finetti's problem the optimality holds only for a subset of \lev processes (see \cite{Loeffen_2008}),  it holds for any spectrally positive demand process in our problem (except that the \lev measure needs to have a light tail).  On the other hand, without the ruin, the value function becomes unbounded (even not Lipschitz continuous), and the verification arguments become significantly more difficult (see the proof of Theorem \ref{theorem_main_no_transaction}).

In our problem, with a fixed cost, two threshold levels  for the optimal $(s, S)$-policy need to be found, which is different from the optimal stopping and de Finetti's problems (with the exceptions of \cite{Bayraktar_2013} and \cite{Loeffen_2009_2}), where only one parameter describes the optimal policy.  In general, choosing two parameters is significantly more difficult.  However the scale function helps one achieve this objective.  Whereas this technique remains to be established, similar arguments can be found in stochastic games \cite{Leung_Yamazaki_2011, Hernandez_Yamazaki_2013} where two parameters characterize the equilibrium between two players.

%
%

The rest of this study is organized as follows.  Section \ref{section_model} presents a mathematical model of the problem.  Section \ref{section_review_levy} reviews spectrally one-sided \lev processes and scale functions.  Section  \ref{section_S_s} presents the computation of  the expected total costs under the $(s,S)$-policy via the scale function.  Section \ref{section_candidate} obtains a candidate policy via the continuous/smooth fit principle and shows its existence. Section \ref{section_verification} verifies its optimality.  Section \ref{section_no_setup_cost} studies the case without a fixed ordering cost.  Section \ref{section_numerics} presents numerical results and  Section \ref{section_conclusion} concludes this study.  

Throughout this study,  $x+ := \lim_{y \downarrow x}$ and $x-  := \lim_{y \uparrow x}$ are used to indicate the right- and left-hand limits, respectively.  Superscripts $x^+ := \max(x, 0)$, $f^+(x) := \max(f(x), 0)$, $x^- := \max(-x, 0)$, and $f^-(x) := \max(-f(x), 0)$ are used to indicate positive and negative parts. Finally, we let $\Delta \xi_t := \xi_t - \xi_{t-}$, for any right-continuous process $\xi$.

\section{Inventory Models} \label{section_model}

Let $(\Omega, \mathcal{F}, \p)$ be a probability space on which a stochastic process $D = \{ D_t; t \geq 0\}$ with $D_0 = 0$, which represents the demand of a single item, is defined.  Under the conditional probability $\p_x$, the initial level of inventory is given by $x \in \R$ (in particular, we let $\mathbb{P} \equiv \mathbb{P}_0$).  Hence, the inventory in the absence of control follows the stochastic process
\begin{align}
X_t := x - D_t, \quad t \geq 0. \label{def_X}
\end{align}
We shall consider the case where $D$ is a spectrally positive \lev process, or equivalently $X$ is a spectrally negative \lev process; we will define these processes formally in the next section. Let $\mathbb{F} := \left\{ \mathcal{F}_t; t \geq 0 \right\}$ be the filtration generated by $X$ (or equivalently by $D$).  

An (ordering) policy $\pi := \left\{ L_t^{\pi}; t \geq 0 \right\}$ is given in the form of an impulse control $(T_1^\pi, u^\pi_1; T_2^\pi, u^\pi_2; \ldots)$ with $L_{0-}^\pi =0$ and $L_t^\pi = \sum_{i: T_i^{\pi} \leq t} u_i^\pi$, $t \geq0$, where $\{ T_i^\pi; i \geq 1 \}$ is an increasing sequence of $\mathbb{F}$-stopping times and $u_i^\pi > 0$,  for $i \geq 1$, is an $\F_{T_i^\pi}$-measurable random variable. Corresponding to every policy $\pi$, the (controlled) inventory process is given by $U^\pi = \{U_t^\pi; t \geq 0 \}$ where $U_{0-}^\pi=0$ and
\begin{align*}
U_t^\pi := X_t + L_t^\pi, \quad t \geq 0.
\end{align*}
We assume that the order quantity is continuous and backorders are allowed. In addition, we assume that there is no lead time, perishability, and lost sales.

The problem is to compute, for a given discount factor $q > 0$,  the total expected costs given by
\begin{align*}
v^\pi (x) &:= \E_x \Big[ \int_0^\infty e^{-qt} f (U_{t}^\pi) \diff t +  \sum_{i=1}^\infty e^{-q T_i^\pi} g(u_i^\pi ) \Big] \\
&= \E_x \Big[ \int_0^\infty e^{-qt} f (U^{\pi}_{t}) \diff t +  \sum_{0 \leq t < \infty} e^{-q t} g(\Delta L_t^{\pi} ) 1_{\{ \Delta L_t^{\pi} > 0 \}}\Big], \quad x \in \R, 
\end{align*}
and to obtain an admissible policy that minimizes it, if such a policy exists.  We shall additionally assume that an admissible policy $\pi$ is such that $\int_0^\infty \exp(-qt) f (U^{\pi}_{t}) \diff t$ and $\sum_{i=1}^\infty \exp(-q T_i^\pi) g(u_i^\pi )$ are both well-defined and finite $\p_x$-a.s.

Here, $f : \R \rightarrow \R$ corresponds to the cost of holding and shortage when $x > 0$ and $x < 0$, respectively.  
Regarding $g$, we assume
\begin{align}
g(y) := C y + K, \quad y > 0, \label{g_def}
\end{align}
for some unit (proportional)  cost of the item $C \in \R$ and fixed ordering cost $K > 0$.  We shall study the case $K=0$ separately in Section \ref{section_no_setup_cost}.  
 As in \cite{Bensoussan_2009, Bensoussan_2005}, we assume the following.
\begin{assump}  \label{assump_f_g}
\begin{enumerate}
\item $f$ is continuous and piecewise continuously differentiable with $f(0) = 0$,  and grows (or decreases) at most polynomially (that is to say, there exist $m, k > 0$ and $N \in \mathbb{N}$ such that $|f(x)| \leq k |x|^N$ for all $x \in \R$ such that $|x| > m$).
\item For some $a \in \R$,
\begin{align}
\tilde{f}(x) &:= f(x) + C  q x, \quad x \in \R, \label{def_f_tilde}
\end{align}
is  decreasing and convex on $(-\infty, a)$ and increasing on $(a,\infty)$. 
\item For some $c_0 > 0$ and $x_0 \geq a$, we have $\tilde{f}'(x) \geq c_0$ for a.e.\ $x \geq x_0$.
\end{enumerate}
\end{assump}
As in \cite{Bensoussan_2009, Bensoussan_2005}, this is a crucial assumption for our analysis, and in particular will be used to verify the existence of the optimal policy. 

Finally, the (optimal) value function is written as
\begin{equation}\label{eq:classical-p}
  v(x):=\inf_{\pi\in \Pi}v^\pi(x), \quad x \in \R,
\end{equation}
where $\Pi$ is the set of all admissible policies.  If the infimum is attained by some admissible policy $\pi^* \in \Pi$, then we call $\pi^*$ an \emph{optimal policy}.

\section{Spectrally Negative \lev Processes and Scale Functions} \label{section_review_levy}

Throughout this paper, we assume that the demand process $D$ is a spectrally positive \lev process.  Equivalently, the process $X$ as in \eqref{def_X} is a spectrally negative \lev process.  By the \emph{L\'evy-Khintchine formula} (see e.g.\ Bertoin \cite{Bertoin_1996}),  any \lev process can be characterized by its Laplace exponent.  Here, Assumption \ref{assump_finiteness_mu} we shall assume below allows us to write the Laplace exponent of $X$ as
\begin{align}
\psi(s)  := \log \E \left[ e^{s X_1} \right] =  \gamma s +\frac{1}{2}\sigma^2 s^2 + \int_{(-\infty,0)} (e^{s z}-1 - s z ) \nu (\diff z), \quad s \geq 0, \label{laplace_spectrally_positive}
\end{align}
where $\nu$ is a \lev measure with the support $(-\infty,0)$ that satisfies the integrability condition $\int_{(-\infty,0)} (1 \wedge |z|^2) \nu(\diff z) < \infty$. 

A \lev process has paths of \emph{bounded variation} a.s.\ or otherwise it has paths of \emph{unbounded variation} a.s.  The former holds if and only if $\sigma = 0$ and $\int_{(-1,0)}|z| \, \nu(\diff z) < \infty$; in this case, the expression \eqref{laplace_spectrally_positive} can be simplified to
\begin{align*}
\psi(s)   =  \delta s + \int_{(-\infty, 0)} (e^{s z}-1 ) \nu (\diff z), \quad s \geq 0,
\end{align*}
with $\delta := \gamma - \int_{(-\infty,0)}z\, \nu(\diff z)$.  We exclude the case in which $X$ is the negative of a subordinator (i.e., $X$ is monotonically decreasing a.s.). This assumption implies that $\delta > 0$ when $X$ is of bounded variation.  

Regarding the \lev measure $\nu$, we make the same assumption as Assumption 3.2 of \cite{Bensoussan_2005} so that $\exp (-\beta X_1)$ has a finite moment for some small $\beta > 0$.
\begin{assump}  \label{assump_finiteness_mu}
We assume that there exists a $\bar{\beta} > 0$ such that
\begin{align*}
\int_{(-\infty, -1]} e^{\bar{\beta} |x|} \nu(\diff x) < \infty.
\end{align*}
This guarantees that 
\begin{align}
\mu := \E [X_1]  = \psi'(0+), \label{drift}
\end{align}
is well-defined and finite.
\end{assump}

For the rest of this section, we briefly review the fluctuation theory of the spectrally negative \lev process and the scale function, which will play significant roles in solving the problem.  Note that, unless otherwise stated, Assumption \ref{assump_finiteness_mu} is not required for the results in the next subsection to hold. 
\subsection{Scale functions}
Fix $q > 0$. For any spectrally negative \lev process, there exists a function called  the  $q$-scale function 
\begin{align*}
W^{(q)}: \R \rightarrow [0,\infty), 
\end{align*}
which is zero on $(-\infty,0)$, continuous and strictly increasing on $[0,\infty)$, and is characterized by the Laplace transform:
\begin{align}
\int_0^\infty e^{-s x} W^{(q)}(x) \diff x = \frac 1
{\psi(s)-q}, \qquad s > \lapinv, \label{scale_function_laplace}
\end{align}
where
\begin{equation}
\lapinv :=\sup\{\lambda \geq 0: \psi(\lambda)=q\}. \notag
\end{equation}
Here, the Laplace exponent $\psi$ in \eqref{laplace_spectrally_positive} is known to be zero at the origin and convex on $[0,\infty)$; therefore $\lapinv$ is well-defined and is strictly positive as $q > 0$.   We also define, for $x \in \R$,
\begin{align*}
\overline{W}^{(q)}(x) &:=  \int_0^x W^{(q)}(y) \diff y, \\
Z^{(q)}(x) &:= 1 + q \overline{W}^{(q)}(x),  \\
\overline{Z}^{(q)}(x) &:= \int_0^x Z^{(q)} (z) \diff z = x + q \int_0^x \int_0^z W^{(q)} (w) \diff w \diff z.
\end{align*}
Because $W^{(q)}(x) = 0$ for $-\infty < x < 0$, we have
\begin{align}
\overline{W}^{(q)}(x) = 0, \quad Z^{(q)}(x) = 1  \quad \textrm{and} \quad \overline{Z}^{(q)}(x) = x, \quad x \leq 0.  \label{z_below_zero}
\end{align}

The most well-known application of the scale function can be found in the two-sided exit problem. Let us define the first down- and up-crossing times, respectively, of $X$ by
\begin{align}
\label{first_passage_time}
\tau_b^- := \inf \left\{ t > 0: X_t < b \right\} \quad \textrm{and} \quad \tau_b^+ := \inf \left\{ t > 0: X_t >  b \right\}, \quad b \in \R.
\end{align}
Then, for any $b > 0$ and $x \leq b$,
\begin{align}
\begin{split}
\E_x \left[ e^{-q \tau_b^+} 1_{\left\{ \tau_b^+ < \tau_0^- \right\}}\right] &= \frac {W^{(q)}(x)}  {W^{(q)}(b)}, \\
 \E_x \left[ e^{-q \tau_0^-} 1_{\left\{ \tau_b^+ > \tau_0^- \right\}}\right] &= Z^{(q)}(x) -  Z^{(q)}(b) \frac {W^{(q)}(x)}  {W^{(q)}(b)}, \\
 \E_x \left[ e^{-q \tau_0^-} \right] &= Z^{(q)}(x) -  \frac q {\Phi(q)} W^{(q)}(x).
\end{split}
 \label{laplace_in_terms_of_z}
\end{align}

The continuity and smoothness of the scale function depend on the path variation of $X$.  First, regarding the behaviors around zero, as in Lemmas 3.1 and 3.2 of \cite{Kuznetsov2013}, 
\begin{align}\label{eq:Wq0}
W^{(q)} (0) &= \left\{ \begin{array}{ll} 0, & \textrm{if $X$ is of unbounded
variation,} \\ \frac 1 {\delta}, & \textrm{if $X$ is of bounded variation,}
\end{array} \right. \\
\label{eq:Wqp0}
W^{(q)'} (0+) &:= \lim_{x \downarrow 0}W^{(q)'} (x) =
\left\{ \begin{array}{ll}  \frac 2 {\sigma^2}, & \textrm{if }\sigma > 0, \\
\infty, & \textrm{if }\sigma = 0 \; \textrm{and} \; \nu(-\infty,0) = \infty, \\
\frac {q + \nu(-\infty, 0)} {\delta^2}, &  \textrm{if }\sigma = 0 \; \textrm{and} \; \nu(-\infty, 0) < \infty.
\end{array} \right.
\end{align}
These properties are particularly useful in applying the continuous/smooth fit principle in stochastic control problems. In this paper, we use this to obtain the candidate thresholds $(s^*, S^*)$ of the optimal  policy; see Section \ref{subsection_cont_smooth_fit} below.

Regarding the smoothness on $\R \backslash\{0\}$, we have the following; see \cite{Chan_2009} for more comprehensive results.
\begin{remark} \label{remark_smoothness}
If $X$ is of unbounded variation or the \lev measure does not have an atom, then it is known that $W^{(q)}$ is $C^1(\R \backslash \{0\})$.  Hence, 
\begin{enumerate}
\item $Z^{(q)}$ is $C^1 (\R \backslash \{0\})$ and $C^0 (\R)$ for the bounded variation case, while it is $C^2(\R \backslash \{0\})$ and $C^1 (\R)$ for the unbounded variation case,
\item $\overline{Z}^{(q)}$ is $C^2(\R \backslash \{0\})$ and $C^1 (\R)$ for the bounded variation case, while it is $C^3(\R \backslash \{0\})$ and $C^2 (\R)$ for the unbounded variation case.
\end{enumerate}
\end{remark}
These smoothness results are important in order to apply the It\^o formula where the  (candidate) value function must be $C^2$ (resp.\ $C^1$) for the case of unbounded (resp.\ bounded) variation.

\begin{figure}[htbp]
\begin{center}
\begin{minipage}{1.0\textwidth}
\centering
\begin{tabular}{c}
 \includegraphics[scale=0.4]{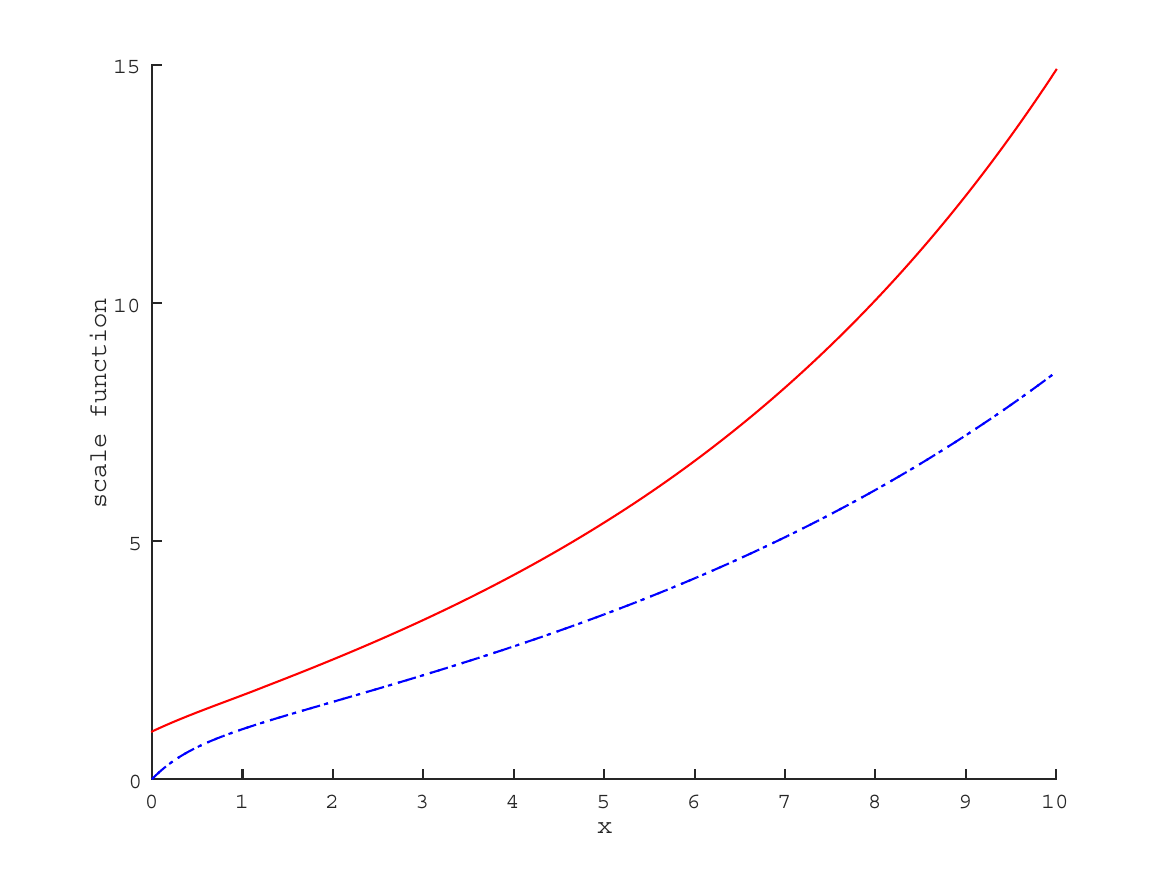}
 \end{tabular}
\end{minipage}
\caption{Plots of the scale function $W^{(q)}$ on $[0, \infty)$. The solid red curve is for the case of bounded variation; the dotted blue curve is for the case of unbounded variation. }  \label{figure_scale_function}
\end{center}
\end{figure}

In Figure \ref{figure_scale_function}, we show  sample plots of the scale function $W^{(q)}$ on $[0, \infty)$. As reviewed in \eqref{eq:Wq0}, the behaviors around zero depend on the path variation of the process. 

\subsubsection{Change of measures}
Fix $\lambda \geq 0$ and define $\psi_\lambda(\cdot)$ as the Laplace exponent of $X$ under $\p^\lambda$ with the change of measure 
\begin{align}
\left. \frac {\diff \p^\lambda} {\diff \p}\right|_{\mathcal{F}_t} = \exp(\lambda X_t - \psi(\lambda) t), \quad t \geq 0;  \label{esscher_transform}
\end{align}
see page 213 of \cite{Kyprianou_2006}.
 Suppose $W_\lambda^{(q)}$ and $Z_\lambda^{(q)}$ are the scale functions associated with $X$ under $\p^\lambda$ (or equivalently with $\psi_\lambda(\cdot)$).  
Then, by Lemma 8.4 of \cite{Kyprianou_2006}, $W_\lambda^{(q-\psi(\lambda))}(x) = \exp(-\lambda x) W^{(q)}(x)$, $x \in \R$,
which is well-defined even for $q \leq \psi(\lambda)$ by Lemmas 8.3 and 8.5 of \cite{Kyprianou_2006}.  In particular, we define
\begin{align}
W_{\Phi(q)}(x) := W_{\Phi(q)}^{(0)}(x) = e^{-\Phi(q) x} W^{(q)}(x), \quad x \in \R, \label{scale_phi_q}
\end{align}
which is known to be an increasing function.

By this change of measure  \eqref{esscher_transform}, one can express the expectation $\E_x \big[ \exp (-q \tau_{0}^-  + v X_{\tau_0^-})\big]$ for $v \geq 0$ using the scale function.  Using this, we can compute the (discounted) moments; for example, by taking the derivative with respect to $v$ and then a limit, we have under \eqref{drift}
\begin{align}
\E_x \big[ e^{-q \tau_{0}^-}  X_{\tau_{0}^-}\big] = \overline{Z}^{(q)}(x) - \mu \overline{W}^{(q)}(x) -  \frac {q - \mu \Phi(q)} {\Phi(q)^2} W^{(q)} (x), \quad x \in \R; \label{first_moment_formula}
\end{align} 
see Proposition 2 of \cite{Avram_et_al_2007} for more detailed computation.

\subsubsection{Martingale properties} \label{subsection_martigale_properties}

 By Proposition 2 of \cite{Avram_et_al_2007} and as in the proof of Theorem 8.10 of \cite{Kyprianou_2006}, the processes 
\begin{align*}
e^{-q (t \wedge \tau^-_{0} \wedge \tau^+_B)} Z^{(q)}(X_{t \wedge \tau^-_{0} \wedge \tau^+_B}) \quad \textrm{and} \quad e^{-q (t \wedge \tau^-_{0} \wedge \tau^+_B)}  R^{(q)} (X_{t \wedge \tau^-_{0} \wedge \tau^+_B}), \quad t \geq 0,
\end{align*}
for any $B > 0$ and $R^{(q)}(y) := \overline{Z}^{(q)}(y) + \mu /q$, $y \in \R$, are $\p_x$-martingales. 

Let $\mathcal{L}$ be the infinitesimal generator associated with
the process $X$ applied to a $C^1$ (resp.\ $C^2$) function $h$ for the case $X$ is of bounded (resp.\ unbounded) variation: for $x \in \R$,
\begin{align} \label{generator}
\begin{split}
\mathcal{L} h(x) &:= \gamma h'(x) + \frac 1 2 \sigma^2 h''(x) + \int_{(-\infty,0)} \left[ h(x+z) - h(x) -  h'(x) z  \right] \nu(\diff z), \\
\textrm{(resp. }\mathcal{L} h(x) &:= \delta h'(x) +  \int_{(-\infty,0)} \left[ h(x+z) - h(x)  \right] \nu(\diff z)\textrm{).}
\end{split}
\end{align}
Thanks to the smoothness of $Z^{(q)}$ and $\overline{Z}^{(q)}$ on $(0,\infty)$ as in Remark \ref{remark_smoothness}, we obtain 
\begin{align}
(\mathcal{L}-q) Z^{(q)}(y) =(\mathcal{L}-q) R^{(q)}(y) = 0, \quad y > 0. \label{martingale_Z_R}
\end{align}
We use these properties in the proof of Lemma \ref{verification_1} below.

\subsubsection{$q$-resolvent (potential) measure} The scale function can express concisely the $q$-resolvent  (potential) density.  As summarized in Theorem 8.7 and Corollaries 8.8 and 8.9 of \cite{Kyprianou_2006} (see also Bertoin \cite{Bertoin_1997}, Emery \cite{Emery_1973}, and Suprun \cite{Suprun_1976}),  we have
\begin{align} \label{resolvent_density}
\begin{split}
\E_x \Big[ \int_0^{\tau_{0}^- \wedge \tau^+_b} e^{-qt} 1_{\left\{ X_t \in \diff y \right\}} \diff t\Big] &= \Big[ \frac {W^{(q)}(x) W^{(q)} (b-y)} {W^{(q)}(b)} -W^{(q)} (x-y) \Big] \diff y,  \quad b > 0, \\
\E_x \Big[ \int_0^{\tau_{0}^-} e^{-qt} 1_{\left\{ X_t \in \diff y \right\}} \diff t\Big] &=\left[ e^{- \Phi(q)y} W^{(q)} (x) -W^{(q)} (x-y) \right] \diff y,  \\
\E_x \Big[ \int_0^\infty e^{-qt} 1_{\left\{ X_t \in  \diff y \right\}} \diff t\Big] &=  \left[  \frac {e^{\Phi(q) (x-y)}} {\psi'(\Phi(q))} -W^{(q)} (x-y) \right] \diff y.
\end{split}
\end{align}
It is clear that these can be used in the computation of inventory costs (see \eqref{lemma_resolvent_eqn}).

The same identities can be obtained when the process $X_t$ is replaced with the 
\emph{running infimum process} $\underline{X}_t := \inf_{0 \leq t' \leq t} X_{ t'}$, $t \geq 0$.  In particular, 
by Corollary 2.2 of \cite{Kuznetsov2013}, for Borel subsets in the nonnegative half-line,
\begin{align*}
\E \Big[ \int_0^{\infty} e^{-qt} 1_{\left\{ - \underline{X}_t \in \diff x \right\}} \diff t \Big] = \frac 1 {\Phi(q)} W^{(q)} (\diff x) -  W^{(q)} (x) \diff x = \frac 1 {\Phi(q)}[\Theta^{(q)}(x) \diff x + W^{(q)}(0) \delta_0(\diff x)],  \end{align*}
where $W^{(q)}(\diff x)$ is the measure such that $W^{(q)}(x) = \int_{[0,x]}W^{(q)}(\diff z)$  (see  \cite[(8.20)]{Kyprianou_2006}) and $\delta_0$ is the Dirac measure at zero.  Here, for all $x > 0$,
\begin{align*} 
\begin{split}
\Theta^{(q)}(x) &:= W^{(q)'} (x+)- \Phi(q) W^{(q)} (x) = e^{\Phi(q)x}W_{\Phi(q)}' (x+) > 0.
\end{split}
\end{align*}
We shall use this function and
\begin{align*}
\overline{\Theta}^{(q)}(x) &:= W^{(q)}(x)    - \Phi(q) \overline{W}^{(q)}(x) > 0,
\end{align*}
  later in the paper. Here the positivity of $\Theta^{(q)}$ holds because $W'_{\Phi(q)}(x+) > 0$ for $x > 0$  and that of $\overline{\Theta}^{(q)}$ holds because it is the integral of $\Theta^{(q)}$. Their positivity will be important in deriving the existence of the optimal solution and the verification of optimality. See the proofs of Proposition \ref{proposition_existence_minimizer} and Lemma \ref{lemma_G_fancy} below.

%
%

\subsection{Examples of scale functions} \label{subsection_examples_scale_function} 
We shall conclude this section with concrete examples of scale functions.  We refer the readers to, e.g., \cite{Hubalek_Kyprianou_2009,Kuznetsov2013} for other examples.

\subsubsection{Brownian motion}  The simplest and nonetheless important example of a \lev process is Brownian motion with a diffusion coefficient $\sigma =1$ and a drift $\delta t$.  In this case, the Laplace exponent \eqref{laplace_spectrally_positive} reduces to $\psi(s) = s^2/2 + \delta s$.  The Laplace transform \eqref{scale_function_laplace} can be analytically inverted; the scale function becomes
\begin{align*}
W^{(q)} (x) = \frac {e^{(\sqrt{\delta^2 + 2 q} - \delta ) x} - e^{-(\sqrt{\delta^2 + 2 q} + \delta ) x}} {\sqrt{\delta^2 + 2 q}}, \quad x \geq 0. 
\end{align*}
\subsubsection{$\alpha$-stable processes} The spectrally negative stable process with index $\alpha \in (1,2)$ has the Laplace exponent $\psi(s) = s^\alpha$.  As in Example 8.2 of \cite{Kyprianou_2006}, its scale function is given by
\begin{align*}
W^{(q)}(x) = \alpha x^{\alpha-1} E'_{\alpha} (q x^\alpha), \quad x \geq 0,
\end{align*}
where $E_\alpha$ is the Mittag-Leffler function of parameter $\alpha$ (which is a generalization of the exponential function).
\subsubsection{Phase-type \lev processes \cite{Asmussen_2004}}  \label{phase_type_case} A spectrally negative \lev process with a finite \lev measure admits a decomposition
\begin{equation}
  X_t  = X_0 + \delta t+\sigma B_t - \sum_{n=1}^{N_t} Z_n, \quad 0\le t <\infty, \label{compound_poisson}
\end{equation}
for some $\delta \in \R$, $\sigma \geq 0$, standard Brownian motion $B=(B_t)_{ t\ge 0}$, a Poisson process $N=(N_t)_{t\ge 0}$ with arrival rate $\rho$ and 
a sequence of i.i.d.\ random variables $Z = (Z_n)_{n = 1,2,\ldots}$. These processes are assumed to be mutually independent.

 The spectrally negative \emph{phase-type} \lev process is a special case where $Z$ is phase-type distributed; a distribution on $(0,\infty)$ is of phase-type  if it is the distribution of the absorption time  in a finite state continuous-time Markov chain consisting of one absorbing state and $m \in\mathbb{N}$ transient states. 
If the phase-type distributed random variable $Z$ is given by a Markov chain with intensity matrix $\bm{T}$ over all  transient states  and the initial distribution $\bm{\alpha}$, then 
the Laplace exponent \eqref{laplace_spectrally_positive} becomes
\begin{align*}
 \psi(s)   = \delta s + \frac 1 2 \sigma^2 s^2 + \rho \left( {\bm \alpha} (s {\bm I} - {\bm{T}})^{-1} {\bm t} -1 \right),
 \end{align*}
which can be extended to $s \in \mathbb{C}$ except at the negative of eigenvalues of ${\bm T}$.

 Suppose $\{ -\xi_{i,q}; i \in \mathcal{I}_q \}$ is the set of the (possibly complex-valued) roots of the equality $\psi(\cdot) = q$ with negative real parts, and if these are assumed distinct, then
the scale function can be written
\begin{align} \label{scale_phase_type}
W^{(q)}(x) = \frac {e^{\Phi(q) x}} {\psi'(\Phi(q))}+ \sum_{i \in \mathcal{I}_q}  \frac 1 {\psi'(-\xi_{i,q})}  e^{-\xi_{i,q}x}, \quad x \geq 0. \end{align}

The set of phase-type \lev processes is dense in the set of all \lev processes, and hence the scale function of any spectrally negative \lev process can be approximated by the scale function of the form  \eqref{scale_phase_type}. See \cite{Egami_Yamazaki_2010_2} for numerical results. 

\subsubsection{Meromorphic \lev processes \cite{Kuznetsov_2010}}
As a variant of the phase-type \lev process, the \emph{meromorphic \lev process} \cite{Kuznetsov_2010} is a type of \lev process whose \lev measure admits a density of the form \begin{align*}
\nu(\diff z) = \sum_{j=1}^\infty p_j \eta_j e^{- \eta_j
|z|} \diff z,  \quad z < 0,
\end{align*}
for some $\{ p_k, \eta_k; k \geq 1 \}$.  Examples include \lev processes in the $\beta$-family as we use for numerical results in Section  \ref{section_numerics} below.
The equation $\psi(\cdot)=q$  has a countable number of negative real-valued roots $\{ -\xi_{k,q}; k \geq 1 \}$ that satisfy the interlacing condition: 
\begin{align*}
\cdots < -\eta_k < -\xi_{k,q} < \cdots < -\eta_2 < -\xi_{2,q} < -\eta_1 < -\xi_{1,q} < 0.
\end{align*}
As discussed in \cite{Kuznetsov2013}, the scale function can be written as 
\begin{align} \label{scale_function_before}
W^{(q)}(x)   = \frac {e^{\Phi(q) x}} {\psi'(\Phi(q))} + \sum_{i=1}^\infty \frac 1 {\psi'(-\xi_{i,q})} e^{-\xi_{i,q} x}, \quad x \geq 0.
\end{align}

 
\section{The $(s,S)$-policy} \label{section_S_s}
In this paper, we aim to prove that the $(s^*,S^*)$-policy is optimal for some $-\infty < s^* < S^* < \infty$.  For arbitrary $-\infty < s < S < \infty$, an $(s,S)$-policy, $\pi_{s,S} := \big\{ L_t^{s,S}; t \geq 0 \big\}$, brings the level of the inventory process $U^{s,S} := X + L^{s,S}$ up to $S$ whenever it goes below $s$.

This process can be defined recursively as follows. First, it moves like the original process $X$ until the first time it goes below $s$:
\begin{align*}
T_1^{s,S} := \tau_s^- = \inf \{ t > 0: X_t < s \}.
\end{align*}
This is immediately pushed up to $S$ (and hence $U^{s,S}_{T_1^{s,S}} = S$) by adding  $u_1^{s,S} := S - X_{T_1^{s,S}}$, and then follows
\begin{align*}
U_t = S + (X_t - X_{T_1^{s,S}}), \quad T_1^{s,S} \leq t < T_2^{s,S}
\end{align*}
where $T_2^{s,S}$ is the first time after $T_1^{s,S}$ the \emph{pre-controlled process}
\begin{align}
\tilde{U}^{s,S}_t := U^{s,S}_{t-} + \Delta X_t, \quad t \geq 0 \label{U_tilde}
\end{align}
 goes below $s$.

The process after $T_2^{s,S}$ can be constructed in the same way. The stopping time $T_i^{s,S}$, for each $i \geq 1$, corresponds to the $i^{th}$ jump of $L^{s,S}$; the $\F_{T_i^{s,S}}$-measurable random variable
\begin{align*} 
u_i^{s,S} := S - \tilde{U}_{T_i^{s,S}}
\end{align*}
 is the corresponding jump size. Clearly, this strategy is admissible: it can be written 
\begin{align*}
L_t^{s,S} = \sum_{i \geq 1: T_i^{s,S} \leq t} u_i^{s,S}, \quad t \geq 0.
\end{align*}

The process $U^{s,S}$ is a strong Markov process.  To see this, the Markov property is clear because, from the construction, the distribution of $U^{s,S}_{t+u}$ only depends on $U^{s,S}_t$ and the increment $X_{t+u}- X_t$, where the latter is independent of $\mathcal{F}_t$. This can be strengthened to the strong Markov property by the right-continuity of the path of $U^{s,S}$ (see, e.g., the proof of Exercise 3.2 of \cite{Kyprianou_2006}, which shows the strong Markov property of a reflected \lev process).

\begin{figure}[htbp]
\begin{center}
\begin{minipage}{1.0\textwidth}
\centering
\begin{tabular}{cc}
 \includegraphics[scale=0.5]{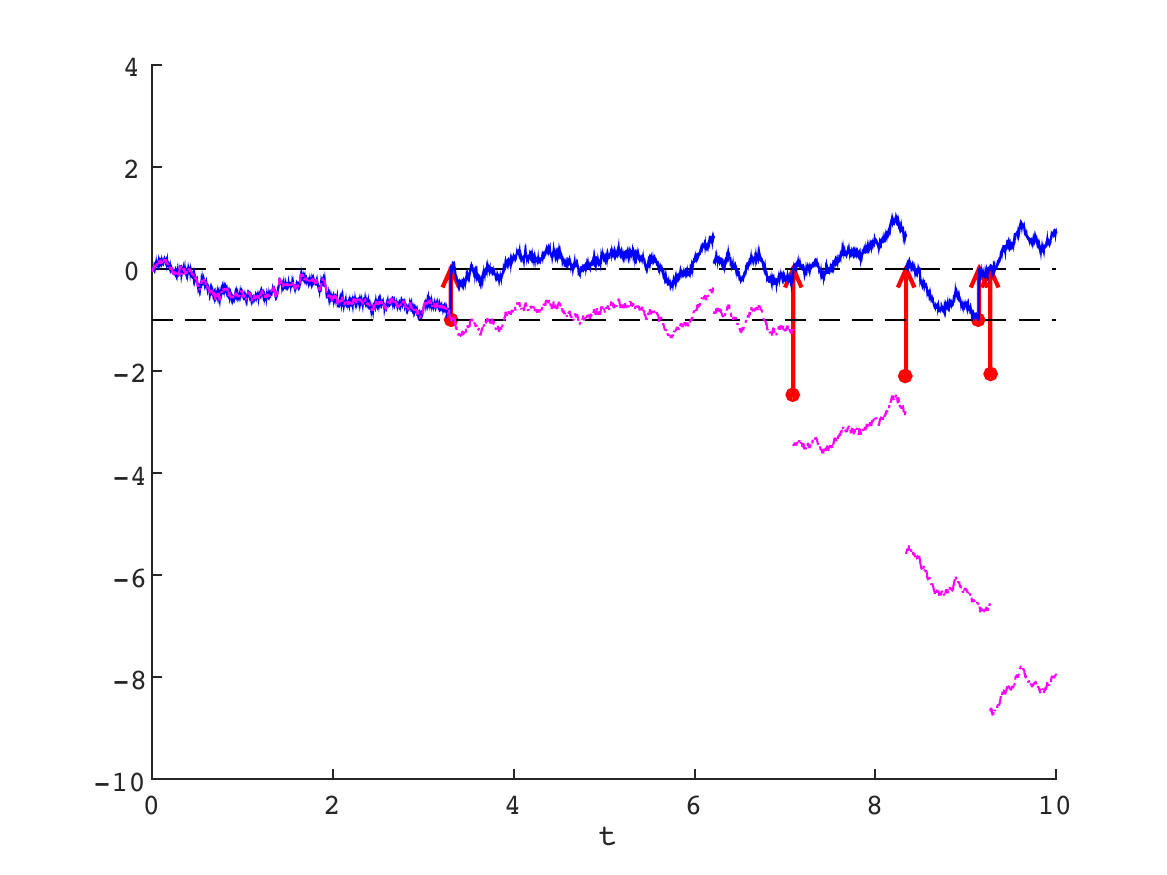} \\ \includegraphics[scale=0.5]{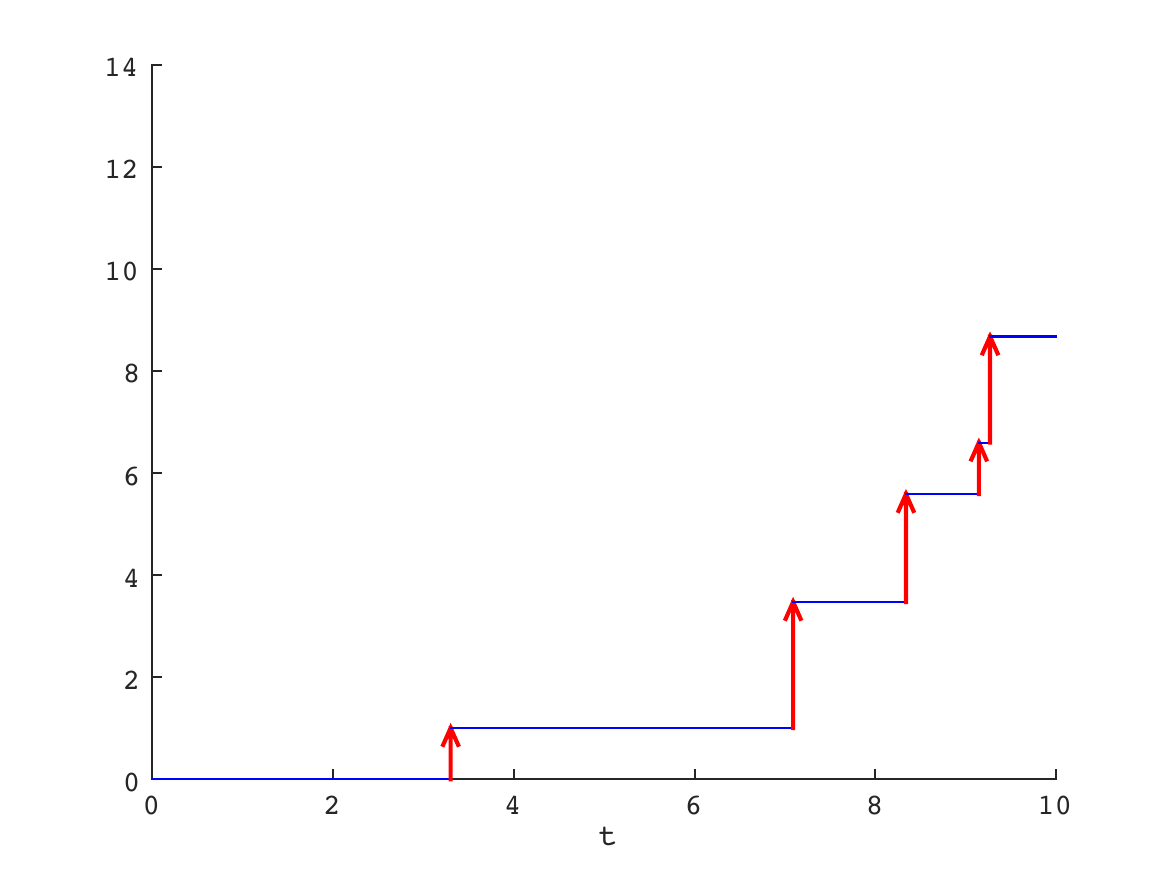}  
\end{tabular}
\end{minipage}
\caption{(Top) sample paths of the underlying process $X$ [pink] and the controlled process $U^{s,S}$ [blue]; (Bottom) the corresponding control process $L^{s,S}$ for $s=-1$, $S=0$.   In the top figure, red arrows show the impulse control that pushes the process up to $S$ whenever the process goes below $s$.   The starting point of the arrow (indicated by a red circle) can be strictly less than $s$ because of the negative jump of the process: it starts at $\tilde{U}^{s,S}_{T_i^{s,S}} := U^{s,S}_{T_i^{s,S}-} + \Delta X_{T_i^{s,S}}$ and ends at $S$. In the bottom figure, the process $L$ becomes the accumulated sum, until $t$, of the increments made by the red arrows.
}  \label{figure_s_S_policy}
\end{center}
\end{figure}

In Figure \ref{figure_s_S_policy}, we show sample paths of the controlled process $U^{s,S}$ and its corresponding control process $L^{s,S}$   for $s=-1$, $S=0$ when the starting value is $x = 0$.  Due to the negative jumps, the process $\tilde{U}^{s,S}$ as in \eqref{U_tilde} can jump to a level strictly below $s$ (and is then immediately pushed to $S$). In the figure, the red arrows show the corresponding impulse control: each time the process goes below $s$, it pushes up to $S$ by adding $u_i^{s,S}$. At the same time, the process $L$ increases by the same amount.

Our objective in this section is to compute the expected total costs under the $(s,S)$-policy, denoted by
\begin{align}
v_{s,S} (x) := \E_x \Big[ \int_0^\infty e^{-qt} f (U^{s,S}_{t}) \diff t +  \sum_{0 \leq t < \infty} e^{-q t} g(\Delta L_t^{s,S} ) 1_{\{ \Delta L_t^{s,S} > 0 \}}\Big], \quad x \in \R. \label{def_v_c}
\end{align}
Toward this end, we compute that of its ``tilted version" (with respect to the unit proportional cost $C$)
\begin{align}
\tilde{v}_{s,S}(x) &:= v_{s,S}(x)  + C x, \quad x \in \R. \label{tilde_relation}
\end{align}
As has been already seen in, e.g., \cite{Bensoussan_2009, Bensoussan_2005}, the computation in the following arguments becomes  simpler if we deal with \eqref{tilde_relation} rather than \eqref{def_v_c}.  The reason will be clear  in later arguments.  In particular, for $x < s$, the slope of $v_{s,S}(x)$ is $-C$ and hence that of $\tilde{v}_{s,S}(x)$ is $0$.  Accordingly, we use $\tilde{f}$ instead of $f$. To see this, note in the verification of optimality that we need to show the positivity of $(\mathcal{L}-q) v_{s,S} (x) + f(x)$, and we have $- q v_{s,S}(x) + f(x) = - q \tilde{v}_{s,S}(x) + \tilde{f}(x)$.

\begin{proposition} \label{proposition_k}
For all $s < S$, \begin{align}
\begin{split}
\tilde{v}_{s,S} (S)   &= \frac  {\Phi(q)} {q  \overline{\Theta}^{(q)}(S-s)}  \Big[ \overline{\Theta}^{(q)}(S-s) \left[ \Psi(s;\tilde{f}) - \frac q {\Phi(q)} \left( K + \frac {C \mu} q \right) \right] + \mathcal{G}(s,S) \Big], \\
\tilde{v}_{s,S} (x) 
&=  \left\{ \begin{array}{ll} - \frac {\overline{\Theta}^{(q)}(x-s)} {\overline{\Theta}^{(q)}(S-s)}\mathcal{G}(s,S)  + \mathcal{G}(s,x)  + \tilde{v}_{s,S}(S), & x \geq s, \\
K +  \tilde{v}_{s,S}(S), & x < s. \end{array} \right. 
\end{split} \label{v_tilde_with_G}
\end{align}
where we define
\begin{align}  
\Psi(s;\tilde{f}) &:=  \int_0^\infty e^{- \Phi(q) y}   \tilde{f}(y+s) \diff y = \int_s^\infty e^{- \Phi(q) (y-s)}   \tilde{f}(y) \diff y, \quad s \in \R, \label{def_psi_phi} \\
\label{def_G}
\mathcal{G}(s,x) 
&:= \Phi(q)   \Psi(s;\tilde{f}) \overline{W}^{(q)}(x-s) + K -  \int_{s}^x W^{(q)} (x-y) \tilde{f}(y) \diff y, \quad x, s \in \R.
\end{align}
\end{proposition}
 
For the rest of the paper, we define $\Psi(s;f)$, $\Psi(s;f')$ and $\Psi(s;\tilde{f}')$ analogously to \eqref{def_psi_phi}. Analytical properties (as in, e.g., Lemma \ref{lemma_bensoussan_result}) of these functions are important in deriving the results of this paper. In particular, for the limiting case $K=0$ as we study in Section \ref{section_no_setup_cost} below, the optimal threshold is given by the root of $\Psi(\cdot;\tilde{f}')$ (see Lemma \ref{lemma_bensoussan_result}).
 
 \begin{remark} \label{remark_finiteness_Psi}
By Assumption \ref{assump_f_g}(1), the functions $\Psi(s;f)$, $\Psi(s;f')$, $\Psi(s;\tilde{f})$ and $\Psi(s;\tilde{f}')$ are finite for any $s \in \R$.
\end{remark}
 

\subsection{Proof of Proposition \ref{proposition_k}}
Recall the down-crossing time $\tau^-_s$ as in \eqref{first_passage_time}. Notice from the construction that,  $\p_x$-a.s.,  $\tilde{U}_t^{s,S} = X_t$ for $0 \leq t \leq \tau_s^-$ and $U_{\tau^-_s}^{s,S} - \tilde{U}_{\tau^-_s}^{s,S} = S- \tilde{U}_{\tau^-_s}^{s,S} = S- X_{\tau^-_s}$ on $\{ \tau^-_s < \infty \}$.
By these, \eqref{g_def}, and  the strong Markov property of  $U^{s,S}$,
the expectation \eqref{def_v_c} must satisfy, for every  $S, x > s$,
\begin{align}
v_{s,S} (x) = \E_x \Big[ \int_0^{\tau_{s}^-} e^{-qt}f(X_{t}) \diff t \Big] + \E_x \left[ e^{-q \tau_{s}^-}  (C (S - X_{\tau_{s}^-})+K) \right] + \E_x \left[ e^{-q \tau_{s}^-}  \right] v_{s,S}(S).   \label{v_recursion}
\end{align}
Define
\begin{align}
k(s,x) &:= \E_{x} \Big[ \int_0^{\tau_{s}^-}  e^{-qt} f(X_{t}) \diff t \Big] -C \E_{x} \left[ e^{-q \tau_{s}^-}  X_{\tau_{s}^-} \right] + K \E_{x} \left[ e^{-q \tau_{s}^-}   \right]  + Cx, \quad x > s. \label{def_k}
\end{align}
Then, using  \eqref{laplace_in_terms_of_z} and \eqref{tilde_relation}, we can write \eqref{v_recursion} as
\begin{align}
\begin{split}\label{v_tilde_2}
\tilde{v}_{s,S} (x) &= \left\{ \begin{array}{ll}k(s,x) + \E_x \big[ e^{-q \tau_{s}^-}  \big] \tilde{v}_{s,S}(S)
= k (s,x) + \left( 1 - \frac q {\Phi(q)} \overline{\Theta}^{(q)}(x-s) \right)\tilde{v}_{s,S} (S),  & \quad x \geq s, \\
 K + \tilde{v}_{s,S} (S), & \quad x < s, \end{array} \right.
\end{split}
\end{align}
where
\begin{align}
\tilde{v}_{s,S} (S) =  \frac  {\Phi(q)} q  \frac {k(s,S) } { \overline{\Theta}^{(q)}(S-s)}, \quad S>s. \label{v_as_f_g}
\end{align}
Here \eqref{v_as_f_g} holds by solving \eqref{v_tilde_2} for $x=S$.  Hence, once we identify the expression for $k(\cdot, \cdot)$, we can compute $\tilde{v}_{s,S}(S)$ and consequently the whole function \eqref{tilde_relation} as well.

The proof consists of evaluating the three expectations in \eqref{def_k}.
The first expectation can be obtained directly by \eqref{resolvent_density}.  We have for $x, s \in \R$,


\begin{align}
\E_{x}\Big[ \int_0^{\tau_{s}^-} e^{-qt} f(X_{t}) \diff t \Big] =
W^{(q)} (x-s)\Psi(s; f) - \int_{s}^x W^{(q)} (x-y) f(y) \diff y,  \label{lemma_resolvent_eqn}
\end{align}
which is well-defined by Remark \ref{remark_finiteness_Psi}.

Here, by the following lemma, we can write \eqref{def_psi_phi} and the integral $\int_{s}^x W^{(q)} (x-y) f(y) \diff y$ interchangeably for $f$ and $\tilde{f}$.
\begin{lemma}  \label{lemma_decomposition_int_w} For $s \in \R$, we have
\begin{align*}
\Psi(s;f) &= \Psi(s;\tilde{f}) - \frac {C q} {\Phi(q)} \left[ \frac 1 {\Phi(q)} + s \right], \\
 \int_{s}^x W^{(q)} (x-y) f(y) \diff y
&= \int_{s}^x W^{(q)} (x-y) \tilde{f}(y) \diff y - C \left[  s Z^{(q)}(x-s) + \overline{Z}^{(q)}(x-s) -x\right], \quad x \in \R.
\end{align*}
\end{lemma}
\begin{proof} See Appendix  \ref{proof_lemma_decomposition_int_w}.
\end{proof}

Using \eqref{lemma_resolvent_eqn} and Lemma \ref{lemma_decomposition_int_w}, we have, for $x,s \in \R$,
\begin{align} \label{running_reward_rewrite}
\begin{split}
\E_{x}\Big[ \int_0^{\tau_{s}^-} e^{-qt} f(X_{t}) \diff t \Big] &= \overline{\Theta}^{(q)}(x-s) \left[ \Psi(s;\tilde{f}) - \frac {C q} {\Phi(q)} \left( \frac 1 {\Phi(q)}+s \right) \right] + C \left[  \overline{Z}^{(q)}(x-s) -(x-s)\right] \\ &+\Phi(q) \overline{W}^{(q)} (x-s) \left[ \Psi(s;\tilde{f}) - \frac {C q} {\Phi(q)^2}  \right]  - \int_{s}^x W^{(q)} (x-y) \tilde{f}(y) \diff y. 
\end{split}
\end{align}

Regarding the second expectation of \eqref{def_k}, combining \eqref{laplace_in_terms_of_z} and \eqref{first_moment_formula} gives the following:
\begin{align} \label{lemma_h_eq}
\E_x \left[ e^{-q \tau_{s}^-}  X_{\tau_{s}^-}\right]
&= \overline{Z}^{(q)}(x-s) - \left(s - \frac \mu q \right) \frac q {\Phi(q)} \overline{\Theta}^{(q)} (x-s) + s - \frac q {\Phi(q)^2} W^{(q)}(x-s), \quad x, s \in \R.
\end{align}

In \eqref{def_k}, substituting \eqref{laplace_in_terms_of_z} (with $x$ replaced with $x-s$),  \eqref{running_reward_rewrite}, and \eqref{lemma_h_eq}, we have the expression
\begin{align*}
k(s,x) 
&= \overline{\Theta}^{(q)}(x-s) \left[ \Psi(s; \tilde{f}) - \frac q {\Phi(q)} \left( K + \frac {C \mu} q \right) \right] + \mathcal{G}(s,x).
\end{align*}
Substituting this in \eqref{v_as_f_g} shows the first equation of \eqref{v_tilde_with_G}. The second equation is immediate by \eqref{v_tilde_2}.

\begin{remark} The same technique can be used to compute the expected costs under the \emph{four parameter band policy} for an extension of the problem with a two-sided control; see the note by Yamazaki \cite{Yamazaki_2014}.
\end{remark}

\section{Candidate policies} \label{section_candidate}
In the previous section, we computed $v_{s,S}$ for arbitrary $s < S$ as in \eqref{v_tilde_with_G}.  Here two different functions for $x > s$ and $x \leq s$ are pasted together at the point $s$, and hence the continuity/smoothness of the function $v_{s,S}$ does not necessarily hold at the point $s$. The principle of smooth/continuous fit chooses the parameters so that the function $v_{s,S}$ becomes continuous/smooth at $s$.  In our problem, as we need to identify two parameters $(s,S)$, we use one additional condition described below.

In this section, we obtain the candidates of $(s,S)$ for the optimal policy.  Toward this end, we choose $(s,S)$ such that  (1) $v_{s,S}(\cdot)$ is continuous (resp.\ differentiable)  at the lower threshold $s$ when $X$ is of  bounded (resp.\ unbounded) variation, and (2) the slope of $v_{s,S}(\cdot)$ at the upper threshold $S$ equals the negative of the unit proportional cost.  We rewrite these two conditions, via the scale function, as $\mathcal{G}(s,S) = 0$ and $\mathcal{H}(s,S) = 0$ where $\mathcal{G}$ is defined as in  \eqref{def_G} and $\mathcal{H}$ is the derivative of $\mathcal{G}$ with respect to the second argument, i.e.,
\begin{align} \label{def_H}
\mathcal{H}(s,x) &:= \frac \partial {\partial x} \mathcal{G}(s,x), \quad x > s.
\end{align} 
We shall then show the existence of the pair $(s^*,S^*)$ that simultaneously satisfy  these two equations.

%

\subsection{Continuous/smooth fit}  \label{subsection_cont_smooth_fit}
We will see that the condition $\mathcal{G}(s,S) = 0$ is equivalent to the so-called \emph{continuous/smooth fit} condition at $s$.   
 It is well-known, in the existence of a diffusion component, that smooth fit holds for impulse control problems; see, e.g., \cite{Guo_Wu_2009}.  On the other hand, for a \lev process of bounded variation, continuous fit may be used alternatively.  This is well-studied particularly for optimal stopping problems; see, e.g., \cite{Egami-Yamazaki-2011,Kyprianou_Surya_2007}.  Here we apply continuous fit for the case $X$ is of bounded variation and smooth fit for the case it is of unbounded variation.
 
 Once $\mathcal{G}(s,S) = 0$ is satisfied, then the second condition $\mathcal{H}(s,S) = 0$ turns out to be equivalent to the condition that the slope of the value function at $S$ equals the negative of the proportional cost $C$, i.e., the slope of $\tilde{v}_{s,S}(\cdot)$ at $S$ is zero.  Note that this condition is commonly used in impulse control to identify optimal $(s,S)$-policies.

For all $x > s$, the function \eqref{def_H} can be written
\begin{align*} 
\mathcal{H}(s,x) &= \Phi(q) \Psi(s; \tilde{f})  W^{(q)}(x-s)  -
 \frac \partial {\partial x}\int_{s}^x W^{(q)} (x-y) \tilde{f}(y) \diff y \\ &= \Psi(s; \tilde{f}')  W^{(q)}(x-s)  - \int_{s}^x W^{(q)} (x-y) \tilde{f}'(y) \diff y,
\end{align*}
where the last equality holds because integration by parts gives
\begin{align}
\int_{s}^x W^{(q)} (x-y) \tilde{f}(y) \diff y  &=  \overline{W}^{(q)}(x-s) \tilde{f}(s) + \int_{s}^{x} \overline{W}^{(q)}(x-y) \tilde{f}'(y) \diff y,  \quad s,x \in \R, \label{tilde_W_integral_form_change}
\end{align}
and 
\begin{align}
\Psi(s;\tilde{f}) =  \frac {\tilde{f}(s) + \Psi(s;\tilde{f}')} {\Phi(q)}, \quad s \in \R. \label{Psi_transformation}
\end{align}

\begin{proposition}\label{lemma_smoothness_s} 
Suppose $(s,S)$ are such that $\mathcal{G}(s,S) =\mathcal{H}(s,S)  = 0$.  Then, 
\begin{enumerate}
\item $\tilde{v}_{s,S}(\cdot)$ is continuous (resp.\ differentiable) at $s$ when $X$ is of bounded (resp.\ unbounded) variation, 
\item $\tilde{v}_{s,S}' (S) = 0$.
\end{enumerate}
\end{proposition}
\begin{proof}
The proof of Proposition \ref{lemma_smoothness_s}  can be carried out by  a straightforward differentiation of the scale function and its asymptotic behavior near zero as in \eqref{eq:Wq0} and \eqref{eq:Wqp0}.

By taking $x \downarrow s$ in the second equality of \eqref{v_tilde_with_G} and because
\begin{align}
\mathcal{G}(s,s) = K > 0, \quad s \in \R, \label{G_s_positive}
\end{align}
we have
\begin{align*}
\tilde{v}_{s,S} (s+) 
&= - \frac {\overline{\Theta}^{(q)}(0+)} {\overline{\Theta}^{(q)}(S-s)}\mathcal{G}(s,S)  +K + \tilde{v}_{s,S}(S) = - \frac {\overline{\Theta}^{(q)}(0+)} {\overline{\Theta}^{(q)}(S-s)}\mathcal{G}(s,S)  + \tilde{v}_{s,S}(s-).
\end{align*}
Note that $\overline{\Theta}^{(q)}(0+) = 0$ if and only if $X$ is of unbounded variation in view of  \eqref{eq:Wq0}. Hence, the continuity at $x = s$ holds if and only if $\mathcal{G}(s,S)=0$ for the case of bounded variation while it holds automatically for the unbounded variation case. 

For the case of unbounded variation, we further pursue the differentiability at $x=s$. 
Differentiating \eqref{v_tilde_with_G},
\begin{align} \label{derivative_with_G}
\tilde{v}_{s,S}' (x) 
&= - \frac {\Theta^{(q)}(x-s)} {\overline{\Theta}^{(q)}(S-s)}\mathcal{G}(s,S)  + \mathcal{H}(s,x), \quad x > s.
\end{align}
Because $\mathcal{H}(s,s+) = 0$ for the case of unbounded variation by \eqref{eq:Wqp0} and $\tilde{v}_{s,S}'(s-) = 0$, we see that the differentiability holds if and only if $\mathcal{G}(s,S)=0$ as well.

We now turn our attention to the slope at $S$.   If we impose $\mathcal{G}(s,S)=0$ in \eqref{derivative_with_G}, we have $\tilde{v}_{s,S}' (S) 
=  \mathcal{H}(s,S)$.
Hence $\tilde{v}_{s,S}' (S) = 0$ if and only if  $\mathcal{H}(s,S) = 0$.  This completes the proof.\end{proof}

\subsection{Existence of $(s^*,S^*)$} \label{section_existence}
We shall now show that there indeed exists a pair $(s^*, S^*)$ such that $\mathcal{G}(s^*,S^*) = \mathcal{H}(s^*,S^*) = 0$; in the next section we show that this requirement is sufficient to prove the optimality of the $(s^*, S^*)$-policy. We refer the readers to the stochastic games \cite{Leung_Yamazaki_2011} and \cite{Hernandez_Yamazaki_2013} where similar arguments are used to identify a pair of parameters that describe the optimal policy.

Recall $c_0$, $x_0$, $a$ as in Assumption \ref{assump_f_g} and that $\Psi(\cdot;\tilde{f}')$ is equivalent to (4.23) of \cite{Bensoussan_2005} (times a positive constant). By Assumption \ref{assump_f_g}, this satisfies the following.   These results are due to Proposition 5.1 of \cite{Bensoussan_2005} and hence the proof is omitted.  Regarding the connection between $a_0$ and the function $\mathcal{G}(\cdot, \cdot)$, see our later discussion in Section \ref{section_no_setup_cost}.
 
\begin{lemma} \label{lemma_bensoussan_result}
\begin{enumerate}
\item There exists a unique number $a_0 < a$ such that $\Psi(a_0; \tilde{f}') = 0$,  $\Psi(x;\tilde{f}') < 0$ if $x < a_0$ and $\Psi(x; \tilde{f}') > 0$ if $x > a_0$.
\item $\Psi'(x; \tilde{f}') > 0$ for $x \leq a$.
\item $\Psi(x; \tilde{f}')  \geq c_0/\Phi(q)$ for $x \geq x_0$.
\end{enumerate}
\end{lemma}

With $a_0$ defined above, we show that the desired pair $(s^*, S^*)$ exists, and in particular $s^*$ lies on the left-hand side and $S^*$ lies on the right-hand side of $a_0$.

\begin{proposition} \label{proposition_existence_minimizer}
There exist $s^* < a_0$ and  $S^* > a_0$ such that  
\begin{align} \label{s_inf_S}
s^* := \max \left\{ s < a_0: \min_{S \geq s} \mathcal{G}(s,S) = 0 \right\} \quad \textrm{and} \quad S^* \in \arg \min_{S \geq s^*} \mathcal{G}(s^*,S),
\end{align}
 with $\mathcal{H}(s^*,S^*) = \mathcal{G}(s^*,S^*) = 0$ where $S^*$ is defined to be any value of $S$ such that $\mathcal{G}(s^*,S)$ is minimized.
\end{proposition}

%

\begin{proof}[Proof of Proposition \ref{proposition_existence_minimizer}] We first rewrite \eqref{def_G} and \eqref{def_H} as integrals of $\Psi(\cdot; \tilde{f}')$ so as to use Lemma \ref{lemma_bensoussan_result} efficiently.
\begin{lemma} \label{lemma_about_G}
We have
\begin{align*}
\mathcal{G}(s,x)  &=  \int_{s}^{x}  \Psi(y; \tilde{f}') \overline{\Theta}^{(q)}(x-y)  \diff y  + K,  \quad x, s , \in \R, \\
\mathcal{H}(s,x)
&= \Psi(x; \tilde{f}')  W^{(q)}(0)  + \int_{s}^x \Psi(y; \tilde{f}')  \Theta^{(q)}(x-y)   \diff y, \quad x > s.
\end{align*}
\end{lemma}
\begin{proof} See Appendix \ref{proof_lemma_about_G}.
\end{proof}

Second, we obtain the asymptotic behavior of $\mathcal{G}$ as follows.
\begin{lemma} \label{lemma_G_fancy}
\begin{enumerate}
\item For every fixed $s \in \R$,  $\lim_{S \uparrow \infty}\mathcal{G}(s,S) = \infty$.  \item For every fixed $S \in \R$,  $\lim_{s \downarrow -\infty}\mathcal{G}(s,S) = -\infty$.
\end{enumerate}
\end{lemma}
\begin{proof} See Appendix \ref{proof_lemma_G_fancy}.
\end{proof}

We are now ready to show the existence of $(s^*,S^*)$.  While this is shown analytically  below, numerical plots of $\mathcal{G}$ and $\mathcal{H}$ in Figure \ref{figure_G_H} of Section \ref{section_numerics} are certainly helpful in understanding the proof.

First we recall the definition of $a_0$ as in Lemma  \ref{lemma_bensoussan_result}(1) and show that we can focus on $s$ smaller than $a_0$. To see this,
for any $S > s \geq a_0$, Lemmas \ref{lemma_bensoussan_result}(1) and \ref{lemma_about_G} and the positivity of $\Theta^{(q)}$ imply that $\mathcal{H}(s,S) > 0$ uniformly; this together with \eqref{G_s_positive} implies that the function $S \mapsto \mathcal{G}(s,S)$ starts at $K > 0$ (at $S=s$) and increases monotonically as $S \uparrow \infty$ while  never touching the x-axis.  

Let us now start at $s = a_0$ and consider decreasing the value of $s$ toward $-\infty$.  Fix $s < a_0$.
By Lemma \ref{lemma_G_fancy}(1), there exists a global minimizer $S(s) \in \arg \min_{S \geq s} \mathcal{G}(s,S)$.  In addition, for any $S \leq a_0$, $\mathcal{H}(s,S) < 0$ and hence $S(s) > a_0$, which also means that a local minimum is attained at $S(s)$ (hence $\mathcal{H}(s, S(s)) = 0$).

Regarding the function $s \mapsto \min_{S \geq s} \mathcal{G}(s,S)$ on $(-\infty, a_0)$, we have the following three properties:
\begin{enumerate}
\item It decreases as $s$ decreases.  Indeed, $\partial \mathcal{G}(s,x) / {\partial s}  =  -\Psi(s; \tilde{f}') \overline{\Theta}^{(q)}(x-s)$,  
which is positive for $s < a_0$ by Lemma \ref{lemma_bensoussan_result}(1). Hence for any $s_1 < s_2 < a_0$, $\min_{S \geq s_1} \mathcal{G}(s_1,S) \leq  \min_{S \geq s_2} \mathcal{G}(s_1,S) < \min_{S \geq s_2} \mathcal{G}(s_2,S)$.  \item It is continuous thanks to the continuity of $\mathcal{G}(s,S)$ with respect to both variables. 
\item Lemma \ref{lemma_G_fancy}(2) implies, for sufficiently small $s$, that $\min_{S \geq s} \mathcal{G}(s,S) < 0$. 
\end{enumerate}
By these three properties, there must exist  $(s^*, S^*)$ such that \eqref{s_inf_S} holds and $\mathcal{G}(s^*,S^*) = \mathcal{H}(s^*, S^*) = 0$.  Note that while $s^*$ is unique by construction, $S^*$ may not be unique in the sense that the infimum of $S \mapsto \mathcal{G}(s^*, S)$ can possibly be attained at multiple values. Moreover, by construction, $s^* < a_0$ and $S^* > a_0$ (because $\mathcal{H}(s^*,S)$ is negative for $S \in (s^*, a_0)$).  This completes the proof of  Proposition \ref{proposition_existence_minimizer}.
\end{proof}

\section{Verification of optimality} \label{section_verification}
With the $(s^*,S^*)$ obtained in Proposition \ref{proposition_existence_minimizer}, we shall show that the $(s^*,S^*)$-policy is optimal.   By substituting  $\mathcal{G}(s^*,S^*)=0$ in \eqref{v_tilde_with_G}, 
\begin{align}
\tilde{v}_{s^*,S^*}(S^*) &= \frac {\Phi(q)} q \Psi(s^*; \tilde{f})  - K  - \frac  {C \mu} q, \label{v_tilde_final_fixed} \\
\tilde{v}_{s^*,S^*} (x) 
&= \mathcal{G}(s^*,x)  + \tilde{v}_{s^*,S^*} (S^*), \quad x \in \R. \label{v_tilde_final}
\end{align}
Substituting this and by the definition of $\mathcal{G}(\cdot,\cdot)$ as in \eqref{def_G},
\begin{align}  \label{v_tilde_optimal_tilde}
\tilde{v}_{s^*,S^*} (x) 
&= \frac {\Phi(q)} q  \Psi(s^*; \tilde{f})Z^{(q)}(x-s^*)   -\int_{s^*}^x W^{(q)} (x-y) \tilde{f}(y) \diff y -  \frac {C \mu} q, \quad x \in \R.
\end{align}
The function $v_{s^*,S^*}(\cdot)$ can be recovered by \eqref{tilde_relation}.
By Lemma \ref{lemma_decomposition_int_w}, we can also write
\begin{align} \label{v_tilde_optimal}
\begin{split}
v_{s^*,S^*} (x) 
&= \left( \frac {\Phi(q)} q  \Psi(s^*; f) + \frac C {\Phi(q)}\right) Z^{(q)}(x-s^*)   - C  \left( \overline{Z}^{(q)}(x-s^*)  + \frac \mu q \right) - \int_{s^*}^x W^{(q)} (x-y) f(y) \diff y.
\end{split}
\end{align}
Notice that  \eqref{v_tilde_final}, \eqref{v_tilde_optimal_tilde} and \eqref{v_tilde_optimal} hold also for $x \leq s^*$ by \eqref{z_below_zero}
 and because $\mathcal{G}(s^*, x) = K$ for $x < s^*$.

We state the main result of the paper for the case $K > 0$.
\begin{theorem} \label{theorem_main}
The $(s^*,S^*)$-policy is optimal over all admissible policies $\Pi$ and the value function \eqref{eq:classical-p} is given by $v_{s^*,S^*} (\cdot)$ as in  \eqref{v_tilde_optimal}.
\end{theorem}

\begin{remark} \label{remark_smoothness_v}
Recall the generator $\mathcal{L}$ of $X$ as in \eqref{generator}. By Remark \ref{remark_smoothness} and Proposition \ref{lemma_smoothness_s}, the function $v_{s^*,S^*}$ is $C^0 (\R)$ and $C^1(\R \backslash \{ s^*\}) $ (resp.\ $C^1(\R)$ and  $C^2(\R \backslash \{ s^*\}) $) when $X$ is of bounded (resp.\ unbounded) variation.  Moreover, the integral part is well-defined and finite by Assumption \ref{assump_finiteness_mu} and because $v_{s^*, S^*}$ is linear below $s^*$.
Hence, $\mathcal{L} v_{s^*,S^*} (\cdot)$ makes sense anywhere on $\R\backslash \{ s^*\}$.  
\end{remark}

\subsection{Proof of Theorem \ref{theorem_main}}
In order to show Theorem \ref{theorem_main}, it suffices to show that the function $v_{s^*, S^*}$ as in \eqref{v_tilde_optimal} satisfies the QVI (quasi-variational inequalities):
\begin{align} \label{QVI_system}
\begin{split}
&(\mathcal{L}-q) v_{s^*,S^*}(x) + f(x)\geq 0,  \quad x \in \R \backslash \{s^*\},\\
&v_{s^*,S^*}(x) \leq K + \inf_{u \geq 0} \left[ Cu + v_{s^*,S^*}(x+u) \right],  \quad x \in \R, \\
&[(\mathcal{L}-q) v_{s^*,S^*}(x) + f(x)] [v_{s^*,S^*}(x) - K - \inf_{u \geq 0} \left[ Cu + v_{s^*,S^*}(x+u) \right]] = 0, \quad x \in \R \backslash \{s^*\};
\end{split}
\end{align}
 see \cite{Bensoussan_Lions_1984, Bensoussan_2005} and also the proof of verification for the no fixed cost case as in Section \ref{section_no_setup_cost}.  Here one needs to apply an appropriate version of the It\^o formula since the value function $v_{s^*, S^*}$ is not smooth enough at $s^*$ to apply the usual version. 
 For the case $X$ is of unbounded variation, because $v_{s^*, S^*}$ is differentiable with absolutely continuous first derivatives, we apply Theorem 71 of Protter \cite{MR2273672}, which is written in terms of the semimartingale local time (for its definition and existence, see page 216 of \cite{MR2273672}). 
   For the case of bounded variation, we apply the Meyer-It\^o formula as in Theorem 70 of \cite{MR2273672}.  In this case, because $U^\pi$ is again of bounded variation, the semimartingale local time process  is identically zero. Hence the extra term due to the discontinuity of  $v_{s^*, S^*}'(s^*)$ vanishes and has no effects after all. 
%

We first show the first part of the QVI, which can be shown easily thanks to the martingale property as reviewed in Section \ref{subsection_martigale_properties} of the scale function and the fact that $v_{s^*,S^*}$ is linear below $s^*$.
\begin{lemma} \label{verification_1} 
\begin{enumerate}
\item $(\mathcal{L}-q) v_{s^*,S^*}(x) + f(x)= 0$ for $ x > s^*$,
\item $(\mathcal{L}-q) v_{s^*,S^*}(x) + f(x) \geq 0$ for $x < s^*$.
\end{enumerate}
\end{lemma}
\begin{proof} 
(1) First, by  \eqref{martingale_Z_R}, 
 $(\mathcal{L}-q) R^{(q)}(y-s^*) = (\mathcal{L}-q) Z^{(q)}(y-s^*) = 0$ for any $y > s^*$. 
On the other hand, as in the proof of Lemma 4.5 of \cite{Egami-Yamazaki-2010-1}, $(\mathcal{L}-q)  \int_{s^*}^x W^{(q)} (x-y) f(y) \diff y = f(x)$.
Hence in view of  \eqref{v_tilde_optimal},  (1) is proved.

 (2)  
Because $\tilde{v}_{s^*, S^*}(x) = K + \tilde{v}_{s^*, S^*} (S^*)$ for $x < s^*$ and by \eqref{Psi_transformation} and \eqref{v_tilde_final_fixed},
\begin{align*}
(\mathcal{L}-q) v_{s^*,S^*}(x) + f(x)  = - q (K + \tilde{v}_{s^*, S^*} (S^*)) - C \mu + Cq x + f(x) = \tilde{f}(x) - \tilde{f}(s^*) - \Psi(s^*; \tilde{f}').
\end{align*}
This is positive by Assumption \ref{assump_f_g}(2) and Lemma \ref{lemma_bensoussan_result}(1), recalling that $x < s^* < a_0 \leq a$.

\end{proof}

The second part of the QVI is given as follows.
\begin{proposition} \label{verification_2}
For every $x \in \R$, we have $v_{s^*,S^*}(x) \leq K + \inf_{u \geq 0} \left[ Cu + v_{s^*,S^*}(x+u) \right]$.  This inequality holds with equality for $x \leq s^*$.
\end{proposition}
It is clear that showing Proposition \ref{verification_2} is equivalent to showing 
\begin{align}
\tilde{v}_{s^*,S^*}(x) \leq K + \inf_{u \geq 0} \tilde{v}_{s^*,S^*}(x+u). \label{qvi_max_tilde}
\end{align}
  Toward this end, we use the following lemma.

\begin{lemma}  \label{lemma_minimum} The following holds true.
\begin{enumerate}
\item $\tilde{v}_{s^*,S^*}(S^*) = \inf_{x \in \R} \tilde{v}_{s^*,S^*}(x)$.
\item $\tilde{v}_{s^*,S^*}(x)$ is decreasing on $[s^*,a_0]$.
\item $\lim_{x \rightarrow \infty}\tilde{v}_{s^*,S^*}(x) = \infty$.
\end{enumerate}
\end{lemma}
\begin{proof} 
By \eqref{v_tilde_final} and because $S^*$ minimizes $\mathcal{G}(s^*,x)$ over $x \in \R$ as in Proposition \ref{proposition_existence_minimizer} (and because $\mathcal{G}(s^*,x)$ is continuous and in particular constant for $x < s^*$), 
the first claim holds. 

The second claim holds because $\tilde{v}_{s^*,S^*}' (x) = \mathcal{H}(s^*, x) < 0$ on $[s^*, a_0]$ in view of Lemmas \ref{lemma_bensoussan_result}(1) and \ref{lemma_about_G}.      The third claim holds by Lemma \ref{lemma_G_fancy}(1).
\end{proof}
For $x \leq s^*$, Proposition \ref{verification_2} holds immediately by Lemma \ref{lemma_minimum}(1) and because $\tilde{v}_{s^*,S^*}(x)  = \tilde{v}_{s^*,S^*}(s^*) = \tilde{v}_{s^*,S^*}(S^*)+K$.  For $s^* \leq x \leq a_0$, it also holds by  Lemma \ref{lemma_minimum}(2).

\begin{proof}[Proof of Proposition \ref{verification_2} for $x > a_0$]
The proof for $x > a_0$ is the most difficult part of the verification.  However, thanks to Lemmas \ref{verification_1}(1) and \ref{lemma_minimum}(3), we can follow the same arguments as the proof of Theorem 1(iii) of Benkherouf and Bensoussan \cite{Bensoussan_2009}, where they show  that \eqref{qvi_max_tilde} holds also for $x > a_0$.

While they assume that the demand process is a mixture of a Wiener process and a compound Poisson process, only a slight (mostly notational) modification is needed to show that it is also valid for a general spectrally positive \lev demand process. Here we illustrate briefly how this is so. 

The proof of Theorem 1 of \cite{Bensoussan_2009}  uses a contradiction argument to show the claim.  The first step is to assume that  \eqref{qvi_max_tilde} does not hold so that $x' := \min \{ x > a_0: \tilde{v}_{s^*, S^*}(x) - \min_{y \geq x} \tilde{v}_{s^*, S^*}(y) > K \}$ is finite, and then construct the points $(x_1, x_2, x_3)$ such that the following properties hold:
\begin{enumerate}
\item $x_1 < x' \leq  x_2 < x_3$ such that  $\tilde{v}_{s^*,S^*}(x)$ [or equivalently $G_s(x)$ in their paper] attains local minima at $x_1$ and $x_3$ and local maximum at $x_2$,
\item $x_3$ is the smallest minimizer of $\tilde{v}_{s^*,S^*}$ over $(x', \infty)$,
\item $x_2$ is the smallest maximizer of $\tilde{v}_{s^*,S^*}$ over $(-\infty, x_3)$,
\item $x_1$ is the smallest minimizer of  $\tilde{v}_{s^*,S^*}$ over $(-\infty, x')$,
\item $K \leq \tilde{v}_{s^*,S^*}(x_2) - \tilde{v}_{s^*,S^*}(x_3)$, and $\tilde{v}_{s^*,S^*}(x) - \tilde{v}_{s^*,S^*}(y) \leq \tilde{v}_{s^*,S^*}(x_2) - \tilde{v}_{s^*,S^*}(x_3)$ for all $x \leq y \leq x_3$.
\end{enumerate}
The only requirements for these are that $\tilde{v}_{s^*, S^*}$ is continuously differentiable on $(a_0, \infty)$ (which holds by Remark \ref{remark_smoothness_v}) and it tends to $\infty$ (which is verified in Lemma \ref{lemma_minimum}).  Hence, we can construct these points with the same properties in exactly the same way.

 Using $(x_1, x_2, x_3)$ constructed above and the IDE 
\begin{align}
(\mathcal{L}-q) \tilde{v}_{s^*,S^*}(x) + \tilde{f}(x) - C \mu = 0 \label{harmonicity_tilde}
\end{align}
 satisfied on $(a_0, \infty)$ (which holds by Lemma  \ref{verification_1}(1)), a contradiction is derived.   In particular, if $a_0 < x_2 \leq a$ (resp.\ if $x_2 > a$), it is derived by comparing \eqref{harmonicity_tilde} for $x=x_1$ and $x=x_2$ (resp.\ $x=x_2$ and $x=x_3$).  Hence, in doing so, the only difference with \cite{Bensoussan_2009} is that the form of the generator $\mathcal{L}$ is more general (and contains the cutoff function inside the integral).  
 
However, local maxima/minima are attained  at $x = x_1, x_2, x_3$, and hence the first derivatives $\partial  \tilde{v}_{s^*,S^*}(x)/{\partial x}$ vanish and the cutoff function in the integral of $\mathcal{L}$ vanishes. In turn, for $n=1,2,3$,
\begin{align*}
\mathcal{L} \tilde{v}_{s^*,S^*}(x_n) &= \frac 1 2 \sigma^2 \tilde{v}_{s^*,S^*}''(x_n) + \int_{(-\infty,0)} \left[ \tilde{v}_{s^*,S^*}(x_n+z) - \tilde{v}_{s^*,S^*}(x_n)  \right] \nu(\diff z).
\end{align*}
Consequently, we end up having the same expression, and the same results hold.

%
\end{proof}

By Lemma \ref{verification_1} and Proposition \ref{verification_2}, we have shown that $v_{s^*,S^*}$ satisfies the QVI as in \eqref{QVI_system}.  This completes the proof of Theorem \ref{theorem_main}.

\section{The Case with No Fixed Ordering Costs} \label{section_no_setup_cost}

In this section, we study a variant of the problem with $K=0$.  Here, we widen the set of admissible policies to accommodate also the processes containing diffuse increments; we consider the set of $\pi := \left\{ L_t^{\pi}; t \geq 0 \right\}$ given by a nondecreasing, right-continuous, and $\mathbb{F}$-adapted process such that $L_{0-}=0$. With $U_t^\pi := X_t + L_t^\pi$, $t \geq 0$, the problem is to compute the total costs:
\begin{align*}
v^\pi (x) := \E_x \Big[ \int_{[0,\infty)} e^{-qt}  \left( f (U_{t}^\pi) \diff t  + C \diff L^\pi_{t}\right) \Big], \quad x \in \R,
\end{align*}
for some $C \in \R$
and to obtain an admissible policy that minimizes it, if such a policy exists.  For $f$, we again impose Assumption \ref{assump_f_g}.  In addition, an admissible policy $\pi$ is such that $\int_0^\infty \exp(-qt) f (U^{\pi}_{t}) \diff t$ and $\int_{[0,\infty)} \exp(-qt)  \diff L^\pi_{t}$ are both well-defined and finite $\p_x$-a.s.

From the results in the previous section, the optimal policy is easily conjectured.  We have seen, for the case $K >0$, that the $(s^*,S^*)$-policy is optimal for some $s^* < a_0 < S^*$. 
In addition, because $(s^*, S^*)$ are such that $\mathcal{G}(s^*,S^*)  =0$ and, for all $s \in \R$, $\mathcal{G}(s,s) = K \xrightarrow{K \downarrow 0} 0$, the distance between $s^*$ and $S^*$ is expected to shrink as $K$ decreases.  Hence it is a reasonable conjecture for the case $K=0$ that a barrier policy with the lower barrier $a_0$ is optimal.  
 We shall show that it is indeed so.

 \begin{figure}[htbp]
\begin{center}
\begin{minipage}{1.0\textwidth}
\centering
\begin{tabular}{c}
 \includegraphics[scale=0.5]{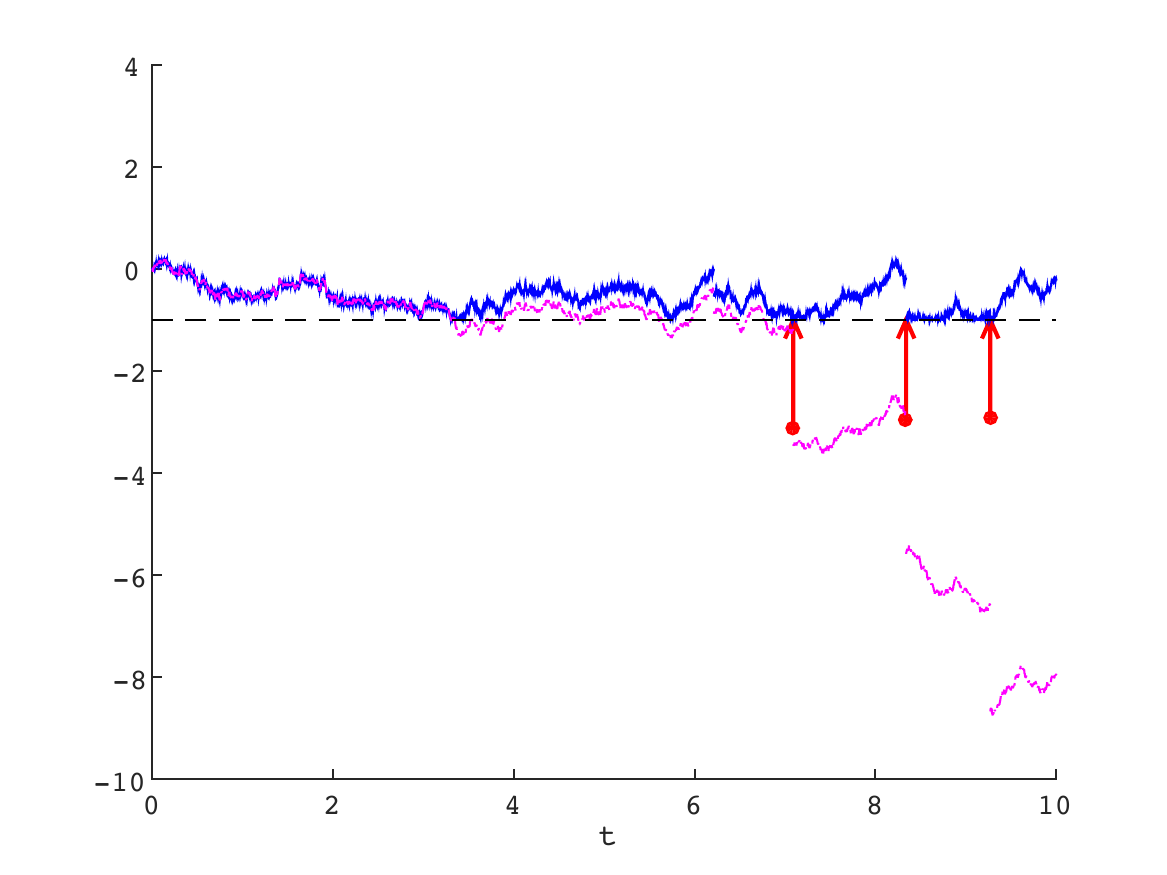} \\  \includegraphics[scale=0.5]{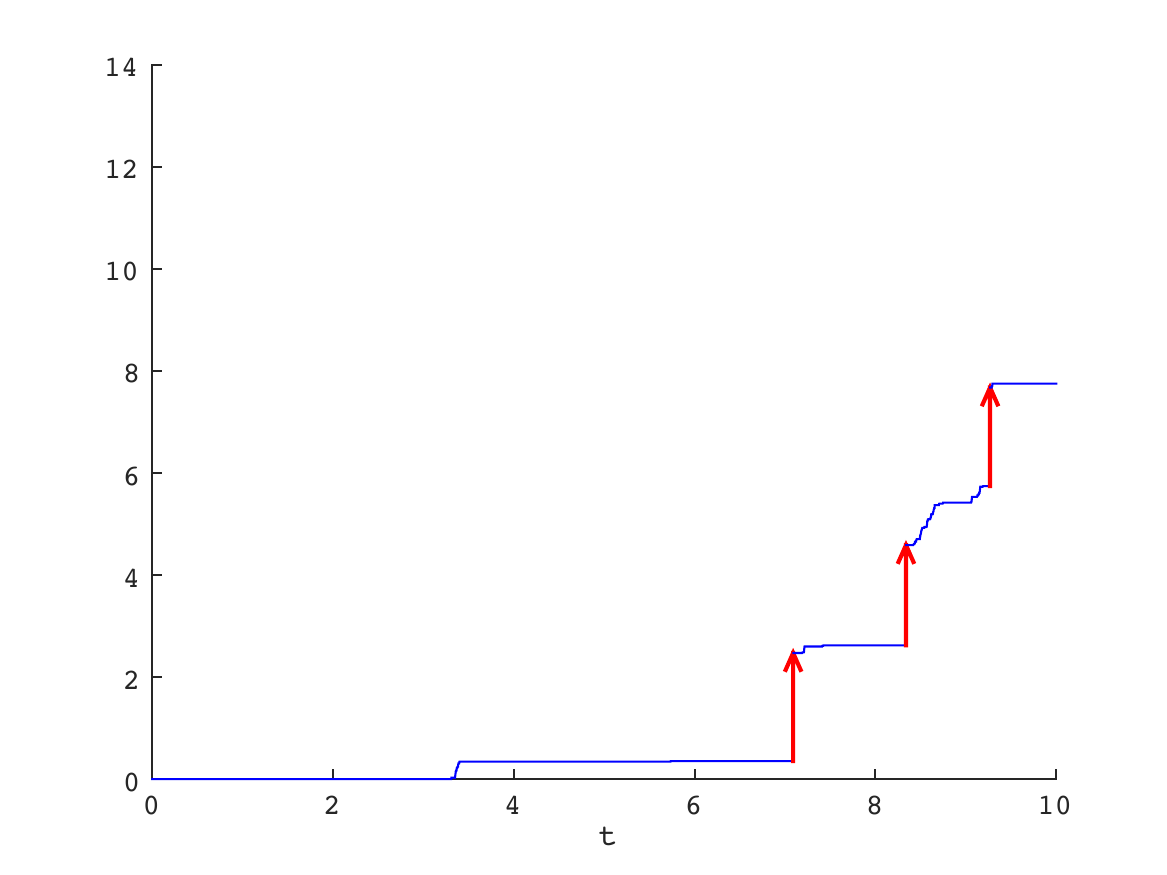}  .
\end{tabular}
\end{minipage}
\caption{(Top) sample paths of the underlying process $X$ [pink] and its reflected process $U^{s}$ [blue]; (Bottom) the corresponding control process $L^{s}$ for $s=-1$.   In the top figure, red arrows show the \emph{discontinuous components} of the singular control that pushes upward the process so that it does not go below $s$.
Each arrow starts at the point $U^{s} + \Delta X$ when it is strictly less than $s$, and ends at $s$. In the bottom figure, the process $L$ becomes the accumulated sum, until $t$, of the increments made by the red arrows \emph{and the continuous increments}.
}  \label{figure_singular_policy}
\end{center}
\end{figure}

Define, for $s \in \R$, 
\begin{align*}
L^{s}_t := \sup_{0 \leq t' \leq t} (s-X_{t'}) \vee 0, \quad t \geq 0.
\end{align*}
The corresponding inventory process $U_t^{s} := X_t + L_t^{s}$ becomes a \emph{reflected \lev process} that always stays at or above $s$.   Figure \ref{figure_singular_policy} shows sample paths of these processes using the same realization of $X$ as in Figure \ref{figure_s_S_policy}; differently from Figure \ref{figure_s_S_policy}, the process $L$ can increase both continuously and discontinuously. 

Our objective is to show that  $v_{a_0}(x) = \inf_{\pi \in \Pi}v^\pi(x)$ for all $x \in \R$ where
\begin{align}
v_{s} (x) := \E_x \left[ \int_{[0,\infty)} e^{-qt}  \left( f (U_{t}^{s}) \diff t  + C \diff L^{s}_{t}\right) \right], \quad x, s\in \R, \label{v_s_no_cost}
\end{align}
and $\Pi$ is the set of all admissible policies.

The fluctuation theory of the reflected \lev process has been well-studied as in, e.g., \cite{Avram_et_al_2007,Pistorius_2004}.  The expression \eqref{v_s_no_cost} can be computed easily via the scale function.
\begin{lemma} \label{lemma_L_singular}
\begin{enumerate}
\item We have $\E_x \big[ \int_{[0,\infty)} \exp(-qt) \diff L^{s}_{t} \big]  = - \overline{Z}^{(q)}(x-s) - \mu / q +  Z^{(q)} (x-s)/\Phi(q)$ for any $x,s \in \R$. 
\item We have $\E_x \big[ \int_0^\infty \exp(-qt)  f (U_{t}^s) \diff t   \big] =  Z^{(q)}(x-s)  {\Phi(q)}  \Psi(s;f) / q -  \int_{s}^x W^{(q)} (x-y) f(y) \diff y$ for any $x,s \in \R$.
\end{enumerate}
\end{lemma}
\begin{proof} See Appendix \ref{proof_lemma_L_singular}.
\end{proof}

Combining the two results above, we can now write
\begin{align} \label{value_function_no_cost_case}
\begin{split}
v_{a_0}(x) 
&=  -  C \left( \overline{Z}^{(q)}(x-a_0) + \frac \mu q  \right) + 
Z^{(q)}(x-a_0) \left( \frac {\Phi(q)} q \Psi(a_0;f) + \frac  C {\Phi(q)} \right)- \int_{a_0}^x W^{(q)} (x-y) f(y) \diff y,
\end{split}
\end{align}
which holds also for $x < a_0$ by \eqref{z_below_zero}.   Here, by Lemma \ref{lemma_decomposition_int_w}, \eqref{Psi_transformation}, and the definition of $a_0$ such that $\Psi(a_0; \tilde{f}')=0$, 
\begin{align*}
\frac {\Phi(q)} q \Psi(a_0;f) + \frac  C {\Phi(q)}  = \frac {f(a_0)} q.
\end{align*}

For the rest of this section, we show the following.
\begin{theorem} \label{theorem_main_no_transaction}
The barrier policy $L^{a_0}$ is optimal and the value function is given by $v_{a_0}(\cdot)$ as in \eqref{value_function_no_cost_case}.
\end{theorem}

\begin{remark}
Recently, Baurdoux and Yamazaki \cite{Baurdoux_Yamazaki_2015} study an extension of this problem where the control is two-sided and in particular show the optimality of a double barrier policy when $f$ is convex.  In another direction, Hern\'andez-Hern\'andez et al.\ \cite{Hernandez_Perez_Yamazaki_2015} consider a version under the condition that $L^\pi$ is absolutely continuous with respect to the \lev measure.
\end{remark}

By Lemma \ref{lemma_decomposition_int_w} and because $\Psi(a_0;\tilde{f}) =  \tilde{f}(a_0) / \Phi(q)$ (by  \eqref{Psi_transformation} and $\Psi(a_0; \tilde{f}')=0$), we can write
\begin{align} \label{value_function_no_cost_case_tilde}
\tilde{v}_{a_0} (x) := v_{a_0} (x) + C x = - \frac {C \mu} q + \frac {Z^{(q)}(x-a_0)} q  \tilde{f}(a_0) - \int_{a_0}^x W^{(q)} (x-y) \tilde{f}(y) \diff y, \quad x \in \R.
\end{align}

%
%

Remark \ref{remark_smoothness} guarantees that $v_{a_0}$ is continuous on $\R$ and  is $C^1 (\R \backslash \{a_0\})$ (resp.\ $C^2 (\R \backslash \{a_0\})$) for the case $X$ is of bounded (resp.\ unbounded) variation.  It turns out that our choice of $s=a_0$ guarantees the smoothness of \eqref{value_function_no_cost_case} (or equivalently \eqref{value_function_no_cost_case_tilde}) at $a_0$ (that is even stronger than the case $K > 0$ as in Section \ref{section_candidate}).  

Because \eqref{tilde_W_integral_form_change} gives $\partial [\int_{a_0}^x W^{(q)} (x-y) \tilde{f}(y) \diff y] /{\partial x}  =  W^{(q)}(x-a_0) \tilde{f}(a_0) + \int_{a_0}^x W^{(q)} (x-y) \tilde{f}'(y) \diff y$, we have
\begin{align}
\tilde{v}_{a_0}' (x) &= {W^{(q)}(x-a_0)}  \tilde{f}(a_0)  - \frac \partial {\partial x}\int_{a_0}^x W^{(q)} (x-y) \tilde{f}(y) \diff y = - \int_{a_0}^x W^{(q)} (x-y) \tilde{f}'(y) \diff y.
\label{v_no_cost_case_derivative}
\end{align}
\begin{lemma} \label{lemma_smoothness_no_cost_case}
The function $v_{a_0}$ is $C^1 (\R)$ (resp.\ $C^2 (\R)$) for the case $X$ is of bounded (resp.\ unbounded) variation.
\end{lemma}
\begin{proof}
By taking $x \downarrow a_0$ in \eqref{v_no_cost_case_derivative}, we see that the differentiability at $a_0$ holds for both bounded and unbounded variation cases.  Furthermore, for the unbounded variation case, by \eqref{eq:Wqp0},
\begin{align*}
\tilde{v}_{a_0}'' (x) 
&=    -\int_{a_0}^x W^{(q)'} (x-y) \tilde{f}'(y) \diff y \xrightarrow{x \downarrow a_0} 0.
\end{align*}
\end{proof}

For the proof of Theorem \ref{theorem_main_no_transaction}, we shall show the following variational inequalities: for all $x \in \R$,
\begin{align} \label{VI_system}
\begin{split}
&(\mathcal{L}-q) v_{a_0}(x) + f(x)\geq 0,  \\
&v_{a_0}' (x)  + C \geq 0,  \\
&[(\mathcal{L}-q) v_{s_0}(x) + f(x)] [v_{a_0}' (x) + C]= 0.
\end{split}
\end{align}
Here, by Lemma \ref{lemma_smoothness_no_cost_case} and because the integral part is well-defined and finite by Assumption \ref{assump_finiteness_mu} and the linearity of $v_{a_0}$ below $a_0$, we confirm that $\mathcal{L} v_{a_0}(x)$ makes sense for all $x \in \R$.

\begin{lemma} \label{verification_1_no_cost} 
\begin{enumerate}
\item $(\mathcal{L}-q) v_{a_0}(x) + f(x)= 0$ for $ x > a_0$,
\item $(\mathcal{L}-q) v_{a_0}(x) + f(x) \geq 0$ for $x \leq a_0$.
\end{enumerate}
\end{lemma}
\begin{proof} 
(1) In view of \eqref{value_function_no_cost_case}, the proof is similar to that of Lemma \ref{verification_1}.

 (2)   By \eqref{value_function_no_cost_case_tilde}, after applying $(\mathcal{L}-q)$ separately to $\tilde{v}_{a_0}$ and $C x$,
\begin{align*}
(\mathcal{L}-q) v_{a_0} (x) + f(x) = -q \Big[ -\frac {C \mu} q + \frac {\tilde{f}(a_0)} q \Big] -  [C\mu - Cq x] + f(x) = \tilde{f}(x) - \tilde{f}(a_0),
\end{align*}
which is negative because $x < a_0 < a$ (by Lemma \ref{lemma_bensoussan_result}(1)) and by Assumption \ref{assump_f_g}(2).
\end{proof}
%

\begin{lemma} \label{verification_2_no_cost} 
We have $\tilde{v}_{a_0}' (x) = 0$ for every $x \leq a_0$ and $\tilde{v}_{a_0}' (x)  \geq 0$ for every $x > a_0$.
\end{lemma}
\begin{proof} 

Recall \eqref{v_no_cost_case_derivative}. 
When $x \leq a_0$, it is clear that $\int_{a_0}^x W^{(q)} (x-y) \tilde{f}'(y) \diff y = 0$.  When $a_0 < x \leq a$, because $\tilde{f}'(y) \leq 0$ for every $a_0 \leq y \leq a$, the claim holds in view of the definition of $\int_{a_0}^x W^{(q)} (x-y) \tilde{f}'(y) \diff y$ as in  \eqref{def_psi_phi}.

Now suppose $x > a$.  Using $W_{\Phi(q)}$ as in  \eqref{scale_phi_q}, we write
\begin{align*}
\begin{split}
\int_{a_0}^x W^{(q)} (x-y) \tilde{f}'(y) \diff y &= e^{\Phi(q)x} \int_{a_0}^x e^{-\Phi(q) y} W_{\Phi(q)} (x-y) \tilde{f}'(y) \diff y \\
&= e^{\Phi(q)x} W_{\Phi(q)} (x-a) \int_{a_0}^x e^{-\Phi(q) y} \tilde{f}'(y) \diff y \\ &+ e^{\Phi(q)x} \int_a^x e^{-\Phi(q) y} [W_{\Phi(q)} (x-y) - W_{\Phi(q)}(x-a)]\tilde{f}'(y) \diff y \\
 &+ e^{\Phi(q)x} \int_{a_0}^a e^{-\Phi(q) y} [W_{\Phi(q)} (x-y) - W_{\Phi(q)}(x-a)] \tilde{f}'(y) \diff y \\
 &\leq e^{\Phi(q)x} W_{\Phi(q)} (x-a) \int_{a_0}^x e^{-\Phi(q) y} \tilde{f}'(y) \diff y,
 \end{split}
\end{align*}
where the last inequality holds because, for a.e.\ $y \in (a,x)$ where $\tilde{f}'(y) \geq 0$, $W_{\Phi(q)} (x-y) - W_{\Phi(q)}(x-a) < 0$ whereas, for a.e.\ $y \in (a_0,a)$ where $\tilde{f}'(y) \leq 0$, $W_{\Phi(q)} (x-y) - W_{\Phi(q)}(x-a) >0$.

Finally, because $\tilde{f}'(y) \geq 0$ for a.e.\ $y \geq x \geq a$,
\begin{align*}
\int_{a_0}^x e^{-\Phi(q) y}\tilde{f}'(y)  \diff y \leq \int_{a_0}^\infty e^{-\Phi(q) y}\tilde{f}'(y)  \diff y = e^{-\Phi(q)a_0}\Psi(a_0; \tilde{f}') = 0. 
\end{align*}
This completes the proof.

\end{proof}

\subsection{Proof of Theorem \ref{theorem_main_no_transaction}}

Lemmas \ref{verification_1_no_cost} and \ref{verification_2_no_cost} show the variational inequalities \eqref{VI_system}.  We shall now give the rest of the proof of Theorem \ref{theorem_main_no_transaction}.

\begin{lemma} \label{lemma_finiteness_resolvent}
For any $x \in \R$,  $v_{a_0} (x) \leq \E_x \left[ \int_0^{\infty} \exp(-qs)   f (X_s) \diff s \right] < \infty$.
\end{lemma}
\begin{proof} See Appendix \ref{proof_lemma_finiteness_resolvent}.
\end{proof}

Now, we are ready to prove Theorem \ref{theorem_main_no_transaction}.  Fix any admissible policy $\pi$ such that $v^\pi(x)$ is finite. Using the standard martingale arguments as in Section 5 of \cite{Oksendal_Sulem_2007} (recall the smoothness of $v_{a_0}$ as in Lemma \ref{lemma_smoothness_no_cost_case}), we obtain
\begin{align}
v_{a_0} (x) \leq \E_x \Big[ \int_{[0, t \wedge T_M^\pi]} e^{-qs}  \left( f (U_{s}^{\pi}) \diff s  + C \diff L^{\pi}_{s}\right) \Big] + \E_x [ e^{-q (t \wedge T_M^\pi)}v_{a_0} (U_{t \wedge T_M^\pi}^\pi)], \quad t, M > 0, \label{v_a_0_bound1}
\end{align}
where we define $T_M^\pi := \inf \{ t \geq 0: |U_t^\pi |> M \}$.

By Lemma \ref{lemma_finiteness_resolvent},  we have that $v_{a_0}(U_{t \wedge T_M^\pi}^\pi) = v_{a_0}(X_{t \wedge T_M^\pi }+ L_{t \wedge T_M^\pi}^\pi)  \leq  \E_{X_{t \wedge T_M^\pi }+ L_{t \wedge T_M^\pi}^\pi} [\int_0^\infty \exp(-q s) f(X_s) \diff s ]$.
Hence, the strong Markov property gives
\begin{align}
&\E_x [e^{-q (t \wedge T_M^\pi )}v_{a_0}(U_{t \wedge T_M^\pi}^\pi)] \leq \E_{x} \Big[\int_{t \wedge T_M^\pi}^\infty e^{-q s} f(L_{t \wedge T_M^\pi}^\pi + X_s) \diff s \Big]. \label{v_a_0_bound2}
\end{align}
Here, for all $s \in [t \wedge T_M^\pi, \infty)$,
\begin{multline*}
 f(L_{t \wedge T_M^\pi}^\pi + X_s) = \tilde{f}(L_{t \wedge T_M^\pi}^\pi + X_s) -  C q  (L_{t \wedge T_M^\pi}^\pi + X_s) \leq \tilde{f} (U_s^\pi)  + \tilde{f} (X_s) -  C q  (L_{t \wedge T_M^\pi}^\pi + X_s)  \\ = f (U_s^\pi) + \tilde{f} (X_s)+   Cq (L_{s}^\pi  - L_{t \wedge T_M^\pi}^\pi),
\end{multline*}
where the first inequality holds because $\widetilde{f}$ is decreasing and increasing by Assumption \ref{assump_f_g}(2) and $X_s \leq L_{t \wedge T_M^\pi}^\pi + X_s \leq U_s^\pi$. Hence, integration by parts gives
\begin{align*} 
\int_{t \wedge T^\pi_M}^\infty e^{-qs}f(L_{t \wedge T_M^\pi}^\pi + X_s) \diff s \leq \int_{t \wedge T^\pi_M}^\infty e^{-qs}f (U_s^\pi) \diff s+ \int_{t \wedge T^\pi_M}^\infty e^{-qs}\tilde{f} (X_s) \diff s +   C \int_{(t \wedge T^\pi_M,\infty)}  e^{-qs}\diff L_s^\pi.
\end{align*}
This together with \eqref{v_a_0_bound1} and \eqref{v_a_0_bound2} gives a bound:
\begin{align*}
v_{a_0}(x) \leq \E_x \Big[ \int_{[0,\infty)} e^{-qs}  \left( f (U_{s}^{\pi}) \diff s  + C \diff L^{\pi}_{s}\right) \Big]  +  \E_x \Big[ \int_{t \wedge T_M^\pi}^\infty e^{-q s} \tilde{f} (X_s) \diff s \Big].
\end{align*}
On the right hand side, the finiteness of $\E_x [\int_0^\infty \exp(-q s) \tilde{f} (X_s) \diff s ] $ can be shown in the same way as the proof  for the finiteness of $\E_x [\int_0^\infty \exp(-q s) f (X_s) \diff s ] $ in Lemma  \ref{lemma_finiteness_resolvent}.
 Hence, by taking $t, M \uparrow \infty$, the claim holds.

\section{Numerical Results} \label{section_numerics}

In this section, we conduct numerical experiments using, for $X$, the spectrally negative \lev process in
 the $\beta$-family introduced by \cite{Kuznetsov_2010_2}. The following definition is due to Definition 4 of \cite{Kuznetsov_2010_2}.
\begin{definition}
A spectrally negative \lev process is said to be in the $\beta$-family if \eqref{laplace_spectrally_positive} is written
\begin{align*}
\psi(z) = \hat{\delta} z + \frac 1 2 \sigma^2 z^2 + \frac \varpi \beta \left\{ B(\alpha + \frac z \beta, 1 - \lambda) - B(\alpha, 1 - \lambda) \right\}
\end{align*}
for some $\hat{\delta} \in \R$, $\alpha > 0$, $\beta > 0$, $\varpi \geq 0$, $\lambda \in (0,3) \backslash \{1,2\}$ and the beta function $B(x,y):=\Gamma(x)\Gamma(y)/\Gamma(x+y)$. 

\end{definition}

The $\beta$-family is a subclass of the meromorphic \lev process and hence the scale function can be computed by the formula \eqref{scale_function_before}.
This process has been receiving much attention recently due to many analytical properties that make many computations possible. In particular, it can approximate tempered stable (or CGMY) processes and is hence suitable to model the price of an asset.  As we discussed in the introduction, the demand is a main determinant of the price; hence it is a reasonable choice for our inventory process $X$.

We suppose $\hat{\delta} = 0.1$, $\lambda = 1.5$, $\alpha=3$, $\beta=1$ and $\varpi = 0.1$.  With this specification, the process has infinitely many jumps in a finite time interval (and has paths of bounded variation), which is not covered in the framework of \cite{Bensoussan_2005}. We consider $\sigma = 0$ and $\sigma = 0.2$ so as to study both the bounded and unbounded variation cases.



We let $q = 0.03$, and, for the inventory cost, we consider the quadratic case $f(x) = x^2$, $x \in \R$.  By straightforward calculation, $a =  - C q/2$, $a_0 = -  {Cq} /2 - {\Phi(q)}^{-1}$,
\begin{align*}
\Psi(s;\tilde{f})  = \frac {2} {\Phi(q)^3} +  \frac {2s + Cq} {\Phi(q)^2} + \frac { s^2+C qs} {\Phi(q)},  \quad   \Psi(s;\tilde{f}') 
= \frac 2 {\Phi(q)^2} + \frac {C q+2s} {\Phi(q)}, \quad s \in \R,
\end{align*}
and for $x \geq s$
\begin{align*}
\int_{s}^x W^{(q)} (x-y) \tilde{f}(y) \diff y
&= \frac {e^{\Phi(q) x}} {\psi'(\Phi(q))} [C q \kappa^{(1)}(s, x; -\Phi(q)) + \kappa^{(2)}(s, x; -\Phi(q)] \\
&-  \sum_{i=1}^\infty B_{i,q} e^{-\xi_{i,q} x} [C q \kappa^{(1)}(s, x; \xi_{i,q}) + \kappa^{(2)}(s, x; \xi_{i,q})], \\
 \int_{s}^x W^{(q)} (x-y) \tilde{f}'(y) \diff y
&= \frac {e^{\Phi(q) x}} {\psi'(\Phi(q))} [C q \kappa^{(0)}(s, x; -\Phi(q)) + 2 \kappa^{(1)}(s, x; -\Phi(q)] \\
&-  \sum_{i=1}^\infty B_{i,q} e^{-\xi_{i,q} x} [C q \kappa^{(0)}(s, x; \xi_{i,q}) + 2 \kappa^{(1)}(s, x; \xi_{i,q})],
\end{align*}
where we define $\kappa^{(n)}(t, t'; \zeta) := \int_t^{t'}  e^{\zeta y} y^n \diff y$ for $t' > t$, $\zeta \in \R$ and $n \geq 0$.

\subsection{Results}  
\begin{figure}[htbp]
\begin{center}
\begin{minipage}{1.0\textwidth}
\centering
\begin{tabular}{cc}
\includegraphics[scale=0.58]{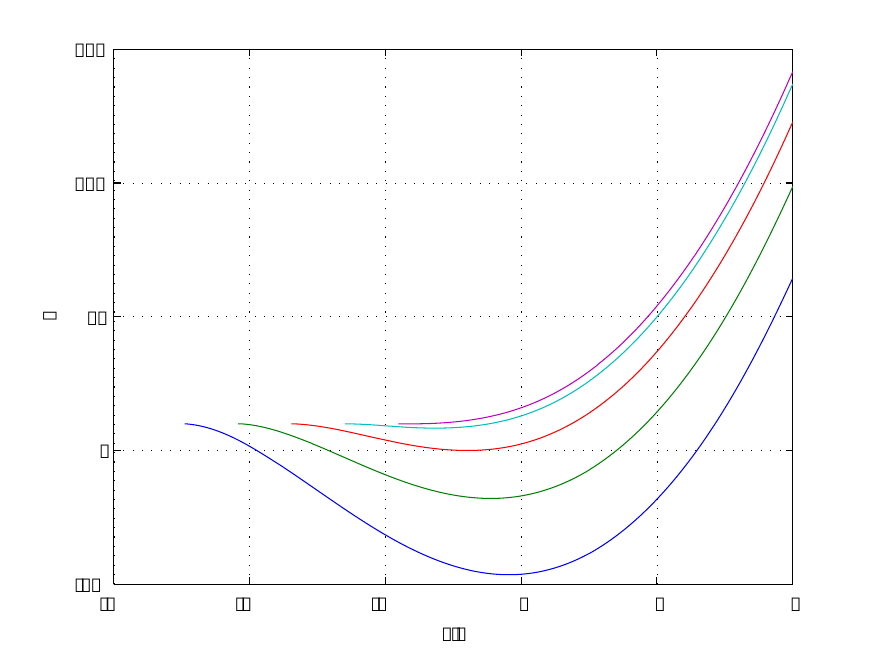}  & \includegraphics[scale=0.58]{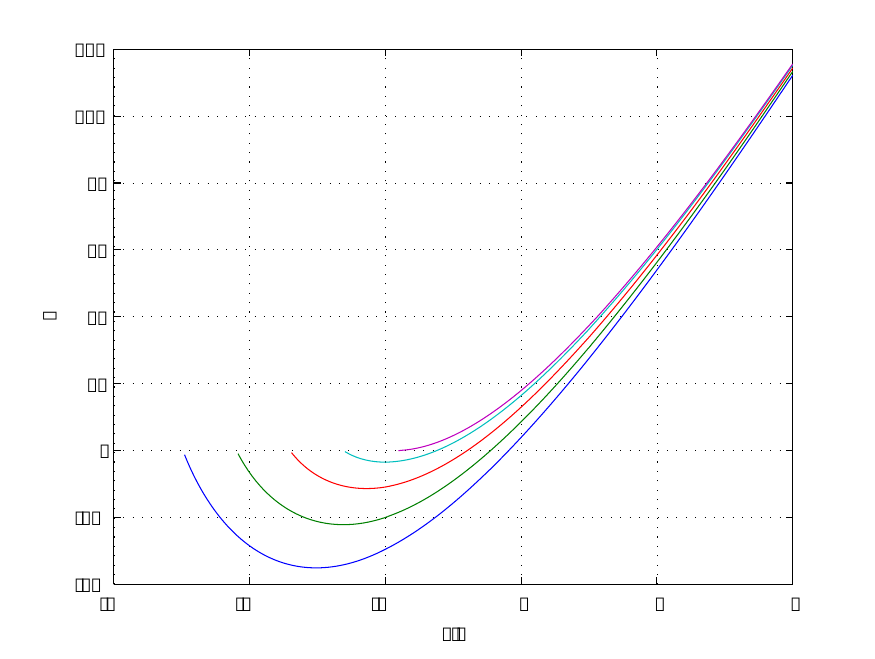} \\
\multicolumn{2}{c}{unbounded variation case ($\sigma > 0$)}\vspace{0.3cm} \\
\includegraphics[scale=0.58]{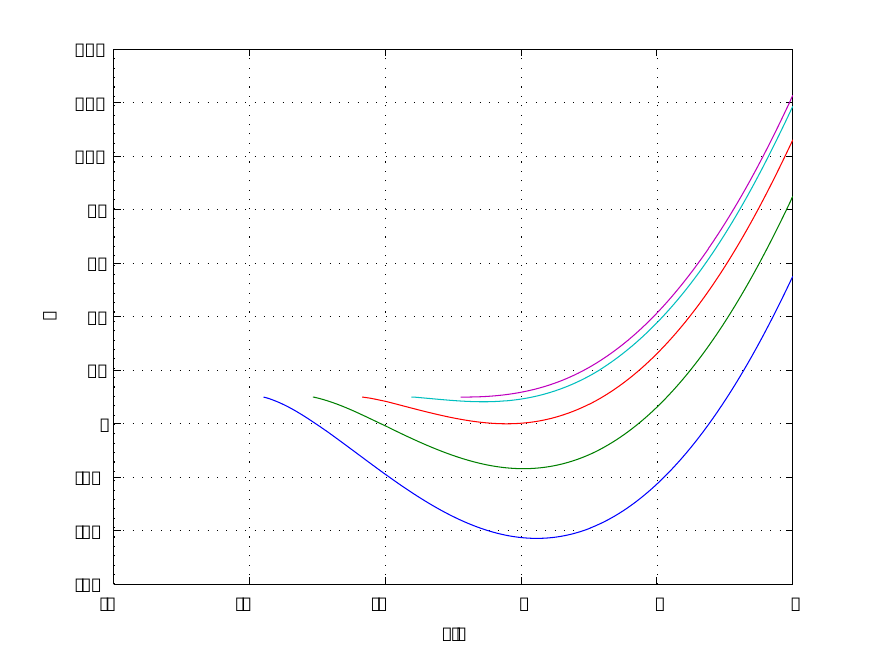}  & \includegraphics[scale=0.58]{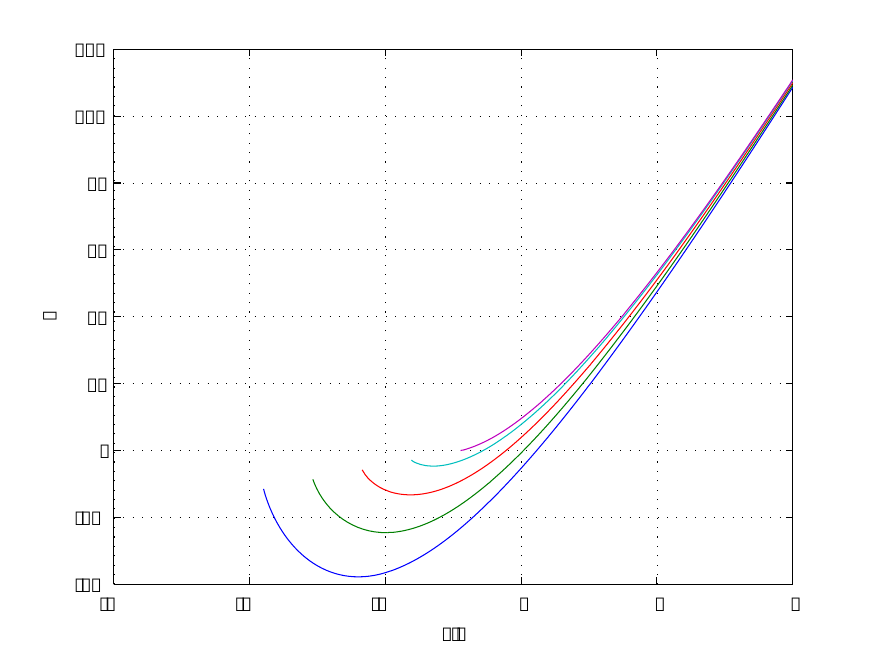} \\
\multicolumn{2}{c}{bounded variation case ($\sigma = 0$)}\vspace{0.3cm} \\
\end{tabular}
\end{minipage}
\caption{Illustrations of $S \mapsto \mathcal{G}(s,S)$ (left) and $S \mapsto \mathcal{H}(s,S)$ (right)  for $s=a_0, (a_0+s^*)/2, s^*, -a_0/2+3s^*/2, -a_0+2s^*$.  Each curve is defined on $[s, \infty)$: it starts at $K > 0$ and eventually goes to $\infty$ as $S \rightarrow \infty$.  When $s = a_0$, it is monotonically increasing.  The value of $\mathcal{G}(s,S)$ is monotone in $s$ and hence as we decrease the value of $s$ from $a_0$, we arrive at a curve that gets tangent to the x-axis (red curve): the starting point becomes $s^*$ and the point touching zero becomes $S^*$.} \label{figure_G_H}
\end{center}
\end{figure}

For the case $K > 0$, the first step is to obtain the pair $(s^*,S^*)$ as in Proposition \ref{proposition_existence_minimizer}.  As is discussed in Section \ref{section_existence}, starting at $s = a_0$ and as we decrease the value of $s$, we arrive at the desired $(s^*,S^*)$ that makes the function $\mathcal{G}(s^*,\cdot)$ tangent to the x-axis at $S^*$.  Figure \ref{figure_G_H} plots $\mathcal{G}(s,\cdot)$ and $\mathcal{H}(s,\cdot)$ for $s=a_0, (a_0+s^*)/2, s^*, -a_0/2+3s^*/2, -a_0+2s^*$.   The lines in red correspond to the desired curve; the starting point becomes $s^*$ and the point touching zero becomes $S^*$. Here we assume $C=K=10$. As it turns out, $ \mathcal{H}(s,\cdot)$ appears to be strictly convex in this example.  Moreover, as is already clear analytically, it starts at zero (for the unbounded variation case) or below zero (for the bounded variation case).  Hence $\mathcal{G}(s,\cdot)$ has a unique global minimum over $S \in [s,\infty)$, which is confirmed to be increasing in $s$ (see the item (1) in the discussion following Lemma \ref{lemma_G_fancy}).  Hence, we apply a bisection method to obtain $s^*$ and $S^*$.

\begin{figure}[htbp]
\begin{center}
\begin{minipage}{1.0\textwidth}
\centering
\begin{tabular}{cc}
 \includegraphics[scale=0.58]{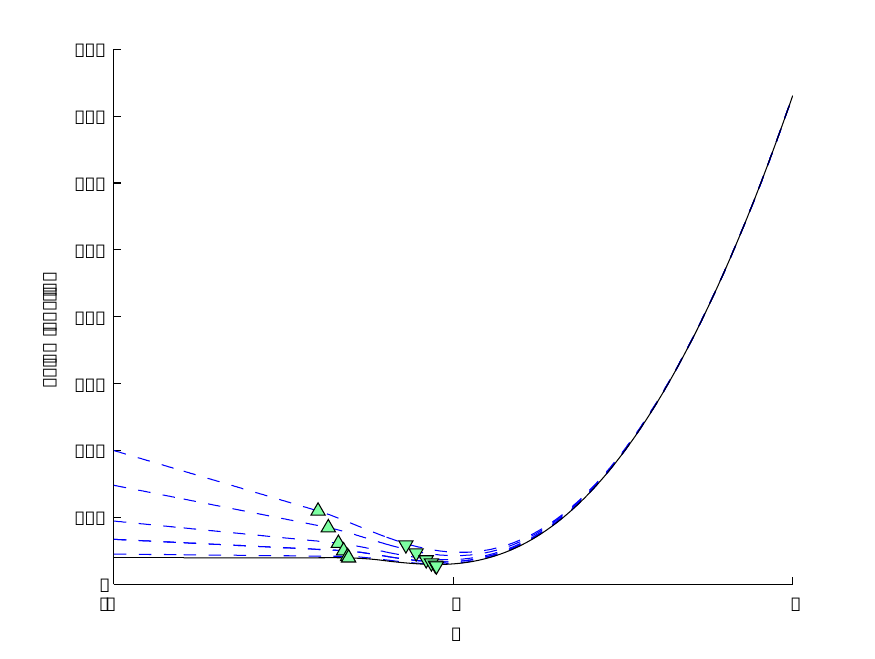}  & \includegraphics[scale=0.58]{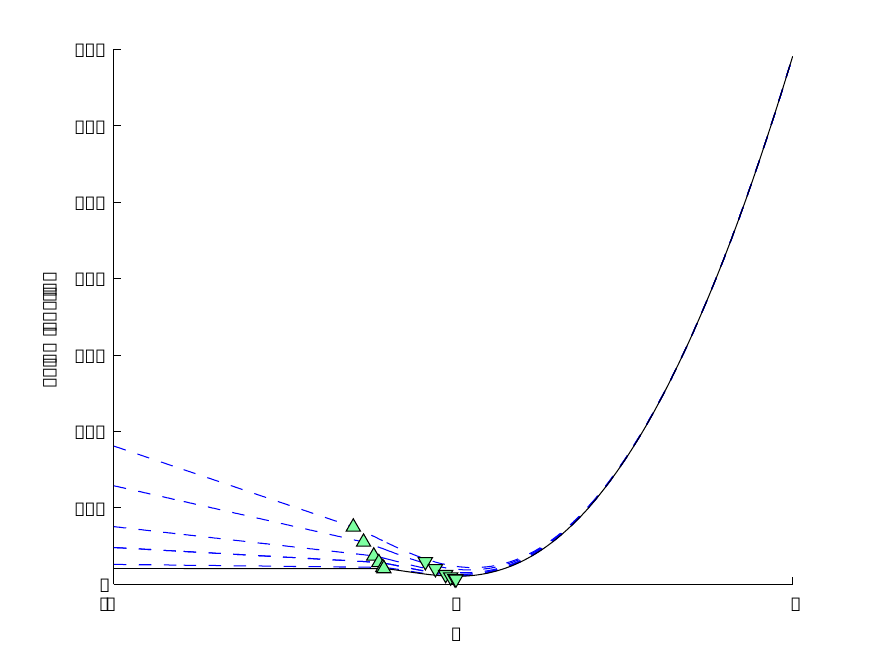}  \\
unbounded variation case ($\sigma > 0$) & bounded variation case ($\sigma = 0$)
\vspace{0.3cm} \\
\end{tabular}
\end{minipage}
\caption{The value functions for various values of the proportional cost $C$.} \label{value_function_proportional}
\end{center}
\end{figure}

With $(s^*,S^*)$ computed instantaneously using the technique addressed above, the value function is computed using \eqref{v_tilde_optimal_tilde}.  In Figure \ref{value_function_proportional}, we first plot it against the initial value $x$  for the unit proportional cost $C=30,20,10,5,1,0$ with the common fixed cost $K=10$.  The triangle signs indicate the points $(s^*, v_{s^*,S^*}(s^*))$ and  $(S^*, v_{s^*,S^*}(S^*))$ for each choice of $C$.  It can be confirmed that the value function is increasing in $C$ uniformly in $x \in \R$.  Moreover, $s^*$  tends to increase as $C$ decreases.  This is consistent with our intuition that one is more eager to replenish as the ordering cost decreases.  We also see  in this plot that $S^*$ also  tends to increase as $C$ decreases.

\begin{figure}[htbp]
\begin{center}
\begin{minipage}{1.0\textwidth}
\centering
\begin{tabular}{cc}
 \includegraphics[scale=0.58]{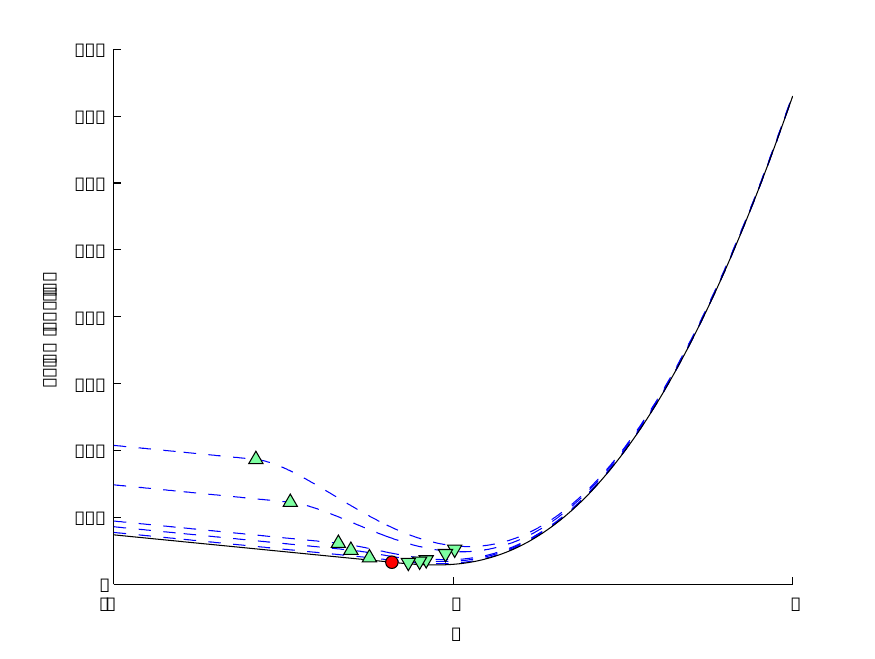} & \includegraphics[scale=0.58]{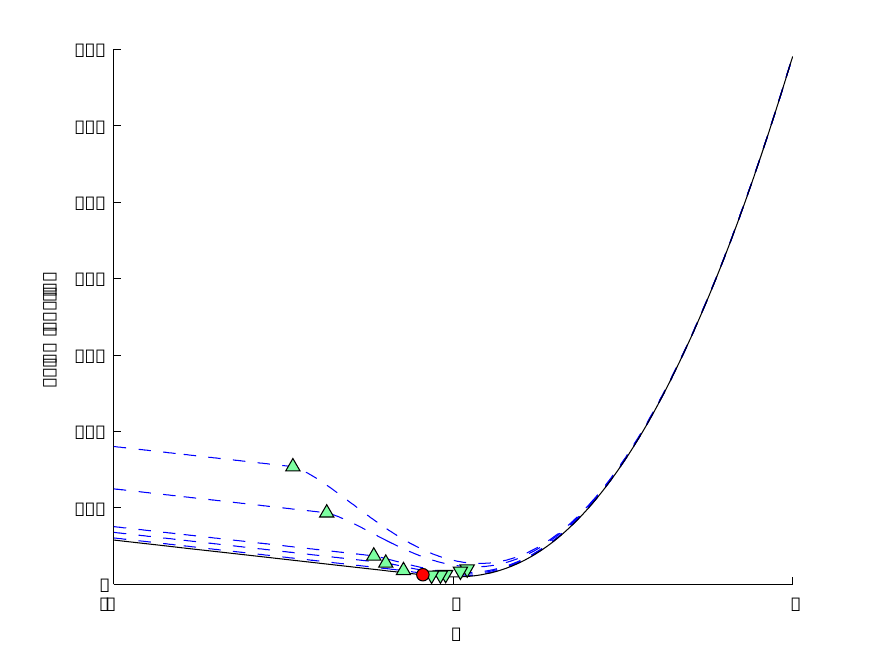}  \\
unbounded variation case ($\sigma > 0$) & bounded variation case ($\sigma = 0$) \vspace{0.3cm} \\
\end{tabular}
\end{minipage}
\caption{The value functions for various values of the fixed cost $K$.} \label{figure_with_diffusion}
\end{center}
\end{figure}

We now consider decreasing the value of the fixed cost $K$ and confirm the convergence to the case $K = 0$ as studied in Section \ref{section_no_setup_cost}.  In Figure \ref{figure_with_diffusion}, we plot the value functions for $K=100,50,10,5,1$ (dotted) along with that for the no-fixed cost case given by \eqref{value_function_no_cost_case} (solid) with the common proportional cost $C=10$.  The circle signs indicate  $(a_0, v_{a_0}(a_0))$. We can confirm that, as the value of $K$ decreases, the value function converges decreasingly to the one for $K=0$.  The convergence of the points $(s^*,S^*)$ to $a_0$ is also confirmed.  Regarding the smoothness of the value function, it appears indeed that it is continuous at $s^*$ for the case of bounded variation while it is differentiable for the case of unbounded variation. Furthermore, the smaller the value of $K$ is, the smoother the value function gets.  This is consistent with Lemma \ref{lemma_smoothness_no_cost_case} where the smoothness holds in a higher order for $K=0$.

\section{Conclusions} \label{section_conclusion}

We have studied the inventory control problem for a general spectrally positive \lev demand process.  We considered both the cases with and without fixed ordering costs. By using the recently developed fluctuation theories of spectrally one-sided \lev processes,  the problem can be solved efficiently and the value function can be written concisely via the scale function. Our numerical experiments show that the computation is fast and accurate; this is a powerful alternative to the existing IDE-based numerical methods, which tend to be difficult when the underlying process has jumps of infinite activity/variation.  

Our approach can potentially be used in other inventory models, and in particular it is of great interest to incorporate, e.g., lead time, perishability  and lost sales.  
The demand for realistic inventory models can potentially contribute to the theory of scale functions.   As insurance problems contributed to the development of the theory of scale functions (see \cite{Albrecher_Renaud_Zhou_2008, Kyprianou_Zhou_2009}), it is expected that pursuing realistic inventory systems will open up new questions on the theory of scale functions and \lev processes.

\appendix

\section{Proofs}

\subsection{Proof of Lemma \ref{lemma_decomposition_int_w}} \label{proof_lemma_decomposition_int_w}
The first claim holds by integration by parts.  For the second claim, we have
\begin{multline*}
\int_{s}^x W^{(q)} (x-y) f(y) \diff y - \int_{s}^x W^{(q)} (x-y) \tilde{f}(y) \diff y=   - C q \int_{s}^{x} y W^{(q)} (x-y) \diff y  \\
=- C q \left[ s \overline{W}^{(q)}(x-s)  +\int_{0}^{x-s} \overline{W}^{(q)} (y) \diff y  \right] =  - C \left[ - x + s Z^{(q)}(x-s) + \overline{Z}^{(q)}(x-s) \right],
\end{multline*}
as desired.

\subsection{Proof of Lemma \ref{lemma_about_G}} \label{proof_lemma_about_G}
We shall show the equation for  $\mathcal{G}(s,x)$ for $x \geq s$ because it is clear that $\mathcal{G}(s,x) = K$ for $x < s$.  The equation for  $\mathcal{H}(s,x)$, $x > s$, is then immediate because it is simply a derivative of $\mathcal{G}(s,x)$.

By integration by parts and \eqref{Psi_transformation},
\begin{align*}
\Psi(s; \tilde{f})  \overline{W}^{(q)}(x-s)    &=  \int_{s}^{x} \left[ \Psi(y; \tilde{f})  W^{(q)}(x-y) - \Psi(y; \tilde{f}')  \overline{W}^{(q)}(x-y)  \right] \diff y \\
&=  \int_{s}^{x} \left[ \frac {\Psi(y;\tilde{f}') } {\Phi(q)}   W^{(q)}(x-y) - \Psi(y; \tilde{f}')  \overline{W}^{(q)}(x-y)  \right] \diff y + \frac {\int_{s}^x W^{(q)} (x-y) \tilde{f}(y) \diff y} {\Phi(q)}  \\
&=\frac 1 {\Phi(q)} \left[ \int_{s}^{x}  \Psi(y; \tilde{f}') \overline{\Theta}^{(q)}(x-y)  \diff y +  \int_{s}^x W^{(q)} (x-y) \tilde{f}(y) \diff y \right].
\end{align*}
Substituting this in \eqref{def_G}, we have the claim.

\subsection{Proof of Lemma \ref{lemma_G_fancy}} \label{proof_lemma_G_fancy}
(1) By Lemmas \ref{lemma_bensoussan_result}(3) and \ref{lemma_about_G}, for any $S \geq s \vee x_0$, 
\begin{align*}
\mathcal{G}(s,S) 
&\geq \frac {c_0} {\Phi(q)}\int_{s \vee x_0}^{S}\overline{\Theta}^{(q)}(S-y)  \diff y  + \int_{s}^{s \vee x_0}  \Psi(y; \tilde{f}') \overline{\Theta}^{(q)}(S-y)  \diff y   + K.
\end{align*}
The claim is now immediate because, for $\epsilon > 0$ and sufficiently large $S$,
\begin{align*}
\int_{s \vee x_0}^{S}\overline{\Theta}^{(q)}(S-y)  \diff y = \int_{0}^{S - s \vee x_0}\overline{\Theta}^{(q)}(y)  \diff y \geq  \int_{\epsilon}^{S - s \vee x_0}\overline{\Theta}^{(q)}(\epsilon)  \diff y \xrightarrow{S \uparrow \infty} \infty,
\end{align*}
where the inequality holds because $\overline{\Theta}^{(q)}$ is positive and monotonically increasing.  

(2) Similarly, if we choose $b < a_0$, we have, for any $s < S \wedge b$,
\begin{align*}
\mathcal{G}(s,S)   =  \int_{s}^{S \wedge b}  \Psi(y; \tilde{f}') \overline{\Theta}^{(q)}(S-y)  \diff y + \int_{S \wedge b}^{S}  \Psi(y; \tilde{f}') \overline{\Theta}^{(q)}(S-y)  \diff y  + K.
\end{align*}
By Lemma \ref{lemma_bensoussan_result}(1) and (2) and because $\Psi(b; \tilde{f}') < 0$, we have, for any $\epsilon > 0$ and sufficiently small $s$,
\begin{multline*}
\int_{s}^{S \wedge b}  \Psi(y; \tilde{f}') \overline{\Theta}^{(q)}(S-y)  \diff y \leq \Psi(b; \tilde{f}') \int_{s}^{S \wedge b}   \overline{\Theta}^{(q)}(S-y)  \diff y \\ =  \Psi(b; \tilde{f}') \int_{S - S \wedge b}^{S-s}   \overline{\Theta}^{(q)}(y)  \diff y \leq \Psi(b; \tilde{f}') \int_{S - S \wedge b + \epsilon}^{S-s}   \overline{\Theta}^{(q)}(S-S \wedge b + \epsilon)  \diff y \xrightarrow{s \downarrow -\infty} -\infty.
\end{multline*}


%


\subsection{Proof of Lemma \ref{lemma_L_singular}} \label{proof_lemma_L_singular}
(1) As in the proof of Theorem 1 of \cite{Avram_et_al_2007}, if we define $\tau_B':= \inf \{ t \geq 0: U_t^0 > B \}$ for $B \in \R$,
\begin{align*}
\E_x \Big[ \int_{[0,\tau_B']} e^{-qt} \diff L^0_{t} \Big] = - l(x) + Z^{(q)}(x) \frac {l(B)} {Z^{(q)}(B)}
\end{align*}
with $l(x) := \overline{Z}^{(q)}(x) + \mu/q - Z^{(q)}(x)/\Phi(q)$.   From Exercise 8.5 of \cite{Kyprianou_2006} and  l'H\^opital's rule,  $\overline{Z}^{(q)}(B)/Z^{(q)}(B) \rightarrow \Phi(q)^{-1}$ as $B \uparrow \infty$.  This together with the monotone convergence theorem gives
\begin{align*}
\E_x \Big[ \int_{[0,\infty)} e^{-qt} \diff L^0_{t} \Big] = - l(x) + Z^{(q)}(x) \lim_{B \uparrow \infty}\frac {l(B)} {Z^{(q)}(B)} 
= - l(x).
\end{align*}
By shifting the initial position of $X$, the proof is complete.

(2) By Theorem 1(i) of \cite{Pistorius_2004}, for every $B > x$ 
\begin{multline*}
\E_x \Big[ \int_0^{\tau'_B} e^{-qt} f^+ (U_{t}^s) \diff t\Big] = \int_s^\infty\Big[ Z^{(q)}(x-s)\frac { W^{(q)} (B-y)} {Z^{(q)}(B-s)} -  W^{(q)} (x-y) \Big]  f^+(y)  \diff y \\
= \int_s^\infty \Big[  Z^{(q)}(x-s) \frac { W^{(q)} (B-s)} {Z^{(q)}(B-s)} e^{-\Phi(q) (y-s)}\frac {W_{\Phi(q)}(B-y)} {W_{\Phi(q)}(B-s)}-  W^{(q)} (x-y) \Big]  f^+(y)  \diff y. 
\end{multline*}
Thanks to Assumption \ref{assump_f_g}(1), the boundedness of $W^{(q)}(\cdot)/Z^{(q)}(\cdot)$ in view of  \eqref{laplace_in_terms_of_z} and because $W_{\Phi(q)}$ is increasing, the integrand of the right hand side is bounded uniformly in $B$ by an integrable function.  Hence, the dominated convergence theorem along with the convergence $W^{(q)}(B)/Z^{(q)}(B) \rightarrow \Phi(q)/q$ (by Exercise 8.5 of \cite{Kyprianou_2006}) yields the result for $f^+$.  The same result holds for $f^-$.  Because both are finite by Assumption \ref{assump_f_g}(1), after summing up these, we have the desired results.


%


\subsection{Proof of Lemma \ref{lemma_finiteness_resolvent}} \label{proof_lemma_finiteness_resolvent}
We shall first prove that 
\begin{align}
\lim_{t, M \uparrow \infty}\E_x [e^{-q (t \wedge \tau_{-M}^- \wedge \tau_{M}^+)} v_{a_0}^+ (X_{t \wedge \tau_{-M}^- \wedge \tau_{M}^+})] = 0 \label{limit_remaining}
\end{align}
  where $\tau^-$ and $\tau^+$ are defined as in \eqref{first_passage_time}.  We have
\begin{align*}
&\lim_{t \uparrow \infty}\E_x \big[e^{-q (t \wedge \tau_{-M}^- \wedge \tau_{M}^+)} v_{a_0}^+ (X_{t \wedge \tau_{-M}^- \wedge \tau_{M}^+}) \big] \\ &= \lim_{t \uparrow \infty}\E_x \big[e^{-q  t} v_{a_0}^+ (X_{t}) 1_{\{ t < \tau_{-M}^- \wedge \tau_{M}^+\}}\big]  +\lim_{t \uparrow \infty}\E_x \big[e^{-q (\tau_{-M}^- \wedge \tau_{M}^+)} v_{a_0}^+ (X_{ \tau_{-M}^- \wedge \tau_{M}^+}) 1_{\{ \tau_{-M}^- \wedge \tau_M^+ \leq t \}}\big]  \\ &= \E_x \big[e^{-q (\tau_{-M}^- \wedge \tau_{M}^+)} v_{a_0}^+ (X_{\tau_{-M}^- \wedge \tau_{M}^+}) 1_{\{ \tau_{-M}^- \wedge \tau_M^+ < \infty\}} \big] \\
&\leq \E_x  \big[e^{-q \tau_{-M}^-} v_{a_0}^+ (X_{\tau_{-M}^-}) 1_{\{ \tau_{-M}^- < \infty\}} \big] + \E_x \big[e^{-q \tau_{M}^+} v_{a_0}^+ (X_{\tau_{M}^+}) 1_{\{ \tau_M^+ < \infty\}} \big],
\end{align*}
where, in the second equality, we used dominated convergence for the limit of the first expectation (because $v_{a_0}^+ (X_{t}) 1_{\{ t < \tau_{-M}^- \wedge \tau_{M}^+ \}}$ is bounded on $[-M, M]$) and monotone convergence for the second expectation.

By the definition of $v_{a_0}$ as in \eqref{v_s_no_cost} and Assumption \ref{assump_f_g}(1), $v_{a_0}^+(z)$  grows at most polynomially as $z \uparrow \infty$ while it is linear below $a_0$.   To see the former, if $U_{t,x}^{a_0}$ is the reflected \lev process that starts at $x \in \R$, it is easy to verify that $U_{t,x + y}^{a_0} \leq U_{t,x}^{a_0} + y$ for all $y > 0$ and $t \geq 0$ a.s.

Hence, for  \eqref{limit_remaining},  it suffices to show 
\begin{align}
\lim_{M \uparrow \infty}\E_x \big[e^{-q \tau_{-M}^-} |X_{\tau_{-M}^-}| 1_{\{ \tau_{-M}^- < \infty\}} \big] = 0 \quad \textrm{and} \quad 
\lim_{M \uparrow \infty}\E_x \big[e^{-q \tau_{M}^+} X_{\tau_{M}^+}^N 1_{\{ \tau_M^+ < \infty\}} \big] = 0, \label{limit_remaining_sides}
\end{align}
for any $N \in \mathbb{N}$.  The latter holds trivially because $X$ does not have positive jumps (and hence $X_{\tau_M^+} = M$ on $\{ \tau_M^+ < \infty\}$ if $x < M$) and $\E_x [\exp(-q \tau_M^+)  1_{\{\tau_M^- < \infty \}}]= \exp(-\Phi(q) (M-x))$, for $M > x$, by Theorem 3.12 of \cite{Kyprianou_2006}. 

By Assumption \ref{assump_finiteness_mu}, we can take a sufficiently small $\beta > 0$ such that $\widetilde{\psi}(\beta) < q$ and  $\int_{(-\infty,-1]} \exp (\beta |x|) \nu (\diff x) < \infty$ where $\widetilde{\psi}(\beta) :=\psi(-\beta)$ is the Laplace exponent of the dual process  $\widetilde{X} := -X$.
Then, as in page 78 of \cite{Kyprianou_2006}, $\{ \exp (\beta \widetilde{X}_t - \widetilde{\psi}(\beta) t), t \geq 0\}$
is a martingale. This means that, for any fixed $M \geq 0$,  $\{ \exp (-\beta X_{t \wedge \tau_{-M}^-} - q (t \wedge \tau_{-M}^-)), t \geq 0 \}$,
is a nonnegative supermartingale.   This together with Fatou's lemma gives, upon taking $t \uparrow \infty$, $\E_x [ \exp (\beta |X_{ \tau_{-M}^-}| - q \tau_{-M}^-) 1_{\{ \tau_{-M}^- < \infty \}}] \leq \exp (-\beta x)$.

Now, for any sufficiently large $M$, we have $(\beta x)^2 \leq \exp (\beta x)$, $x \geq M$, and hence
\begin{align*}
e^{-\beta x} \geq  \beta^2 \E_x \left[e^{- q \tau_{-M}^-}  |X_{ \tau_{-M}^-}|^2   1_{\{ \tau_{-M}^- < \infty \}} \right] \geq \beta^2 M \E_x \left[  e^{- q \tau_{-M}^-} |X_{ \tau_{-M}^-}|  1_{\{ \tau_{-M}^- < \infty \}} \right] \geq 0.
\end{align*}
By taking $M \uparrow \infty$, we see that the first claim of \eqref{limit_remaining_sides} holds and consequently \eqref{limit_remaining} holds.

By the It\^o formula (thanks to the smoothness of $v_{a_0}$ by Lemma \ref{lemma_smoothness_no_cost_case}), Lemma \ref{verification_1_no_cost}  and \eqref{limit_remaining}, 
\begin{align*}
\begin{split}
v_{a_0} (x) &\leq \E_x \Big[ \int_0^{t \wedge \tau^-_{-M} \wedge \tau^+_M} e^{-qs}   f (X_s) \diff s \Big] + \E_x \big[ e^{-q (t \wedge \tau^-_{-M} \wedge \tau^+_M)}v_{a_0}^+ (X_{t \wedge \tau^-_{-M} \wedge \tau^+_M}) \big]  \\ &\xrightarrow{t, M \uparrow \infty} \E_x \left[ \int_0^{\infty} e^{-qs}   f (X_s) \diff s \right].
\end{split}
\end{align*}

Here, we notice that this limit is well-defined and finite.  Indeed, for any $N \in\mathbb{N}$, by the arguments similar to the above, there exist a sufficiently small $\beta' > 0$ (such that $q > \psi(-\beta')$) and large $\underline{M} > 0$,  such that,  for any $t \geq 0$, $\exp [-(q-\psi(-\beta')) t] \exp(-\beta' x) =  \E_x [\exp (-\beta' X_t - qt) ] \geq \E_x [\exp (-\beta' X_t - qt) 1_{\{ X_t < -\underline{M} \}}  ] \geq  (\beta')^N \E_x \left[\exp(- q t)  |X_t|^N 1_{\{ X_t < -\underline{M} \}} \right]$.  Hence, 
\begin{align*}
\E_x \left[ \int_0^{\infty} e^{-qs} |X_s|^N 1_{\{ X_s < - \underline{M} \}} \diff s\right] =   \int_0^{\infty} \E_x \left[  e^{-qs} |X_s|^N 1_{\{ X_s < - \underline{M} \}} \right] \diff s < \infty.
\end{align*}
Similarly, we can show that there exists some $\overline{M} > 0$ such that  $\E_x [ \int_0^{\infty} \exp(-qs) X_s^N 1_{\{ X_s >  \overline{M} \}} \diff s]  < \infty$.  By Assumption \ref{assump_f_g}(1), $\E_x \left[ \int_0^{\infty} \exp(-qs)   f (X_s) \diff s \right]$ is indeed well-defined and finite.

\section*{Acknowledgements}
The author thanks Alain Bensoussan as well as the anonymous referees for their  insightful comments that help improve the presentation of this paper.
K.\ Yamazaki is in part supported by MEXT KAKENHI grant numbers  22710143 and 26800092, JSPS KAKENHI grant number 23310103, the Inamori foundation research grant, and the Kansai University subsidy for supporting young scholars 2014.
\bibliographystyle{abbrv}
\bibliography{dual_model_bib}

\end{document}